\documentclass[letterpaper,12pt,peerreviewca,draftcls]{IEEEtran}
\usepackage{csm16}
\usepackage[margin=1in]{geometry}
\usepackage[nolists,nomarkers,tablesfirst]{endfloat} 
\usepackage{amsmath} 
\usepackage{amsfonts}
\usepackage{circuitikz}
\usetikzlibrary{backgrounds}
\usetikzlibrary{shapes.symbols}
\usetikzlibrary{arrows.meta}
\usepackage{mathtools}
\usepackage{psfrag}

\usepackage{url}
\usepackage{graphicx,xcolor}
\usepackage{verbatim}
\makeatletter
\newcommand{\verbatimfont}[1]{\def\verbatim@font{#1}}%
\makeatother
\verbatimfont{\ttfamily\small}

\newcommand{\bi}{\begin{itemize}}\newcommand{\ei}{\end{itemize}}
\newcommand{\be}{\begin{equation}}\newcommand{\ee}{\end{equation}}
\newcommand{\bee}{\begin{enumerate}}\newcommand{\eee}{\end{enumerate}}
\newcommand{\bea}{\begin{eqnarray}}\newcommand{\eea}{\end{eqnarray}}
\newcommand{\beas}{\begin{eqnarray*}}\newcommand{\eeas}{\end{eqnarray*}}
\newcommand{\bc}{\begin{center}}\newcommand{\ec}{\end{center}}
\usepackage[left,pagewise]{lineno} 
\usepackage[english]{babel} 
\usepackage{blindtext}
\usepackage[compress]{cite}

\title{Dissipativity, reciprocity and passive network synthesis\\
\Large From Jan Willems' seminal Dissipative Dynamical Systems papers to the present day\footnote{\scriptsize This work has been submitted to the IEEE for possible publication. Copyright may be transferred without notice, after which this version may no longer be accessible}}
\author{Timothy H.\ Hughes and Edward H.\ Branford\\
	POC: T.H.Hughes@exeter.ac.uk\\ \today }

\newif\ifPDF \ifx\pdfoutput\undefined\PDFfalse \else\ifnum\pdfoutput > 0\PDFtrue \else\PDFfalse \fi \fi
\ifPDF 
\usepackage[pdftex, plainpages = false, colorlinks=true, linkcolor=black, citecolor = green!50!blue, urlcolor = blue, filecolor=black, pagebackref=false, hypertexnames=false,  pdfpagelabels, draft]{hyperref}
\fi

\begin{document}
\ctikzset{bipoles/resistor/height=.23}
\ctikzset{bipoles/resistor/width/.initial=.6}
\ctikzset{bipoles/capacitor/height/.initial=.5}
\ctikzset{bipoles/capacitor/width/.initial=.13}
\ctikzset{bipoles/length=1.65cm}
\ctikzset{bipoles/cuteinductor/coils=4}
\ctikzset{bipoles/cuteinductor/lower coil height=.15}
\ctikzset{bipoles/cuteinductor/width=.65}
\ctikzset{bipoles/cuteinductor/height=.25}
\ctikzset{bipoles/cuteinductor/coil aspect=0.25}
\maketitle
\CSMsetup

The twin papers Dissipative Dynamical Systems Parts I and II of Jan Willems \cite{JWDSP1, JWDSP2} are considered by many to be the founding papers of the theory of dissipativity within the field of mathematical systems and control. These papers unified earlier work by Kalman, Youla, Anderson, Yakubovich, Popov, and others, relating to electrical network synthesis, optimal control, and stability of interconnected dynamical systems, thereby establishing dissipativity as a standalone concept at the intersection of physics, systems theory, and control engineering. The papers have proved a source of inspiration to many researchers and led to a significant subsequent body of literature, on the dissipativity concept itself (see the sidebar titled ``\nameref{sb:cdac}''), and other significant trends and developments in the field such as port Hamiltonian systems (see \cite{vdsph} and the references therein) and H-infinity optimal control (see the sidebar titled ``\nameref{sb:hinf}''). Indeed, Jan Willems also continued to revisit and to seek to improve upon the dissipativity concept throughout his life, even in the context of finite dimensional linear time-invariant systems. This included an emphasis on addressing the use of an input-state-output formalism and assumptions of controllability and observability in the founding papers (see, for instance, \cite{JWDDS}). This relates closely to the emphasis in his behavioral theory of moving away from an input-output perspective as the starting point for system modelling and analysis, and the establishment of system theoretic properties such as controllability, stabilisability, and the state of a system from a trajectory-level (representation-free) perspective as opposed to properties of a specific state-space representation of the system. Examples abound of systems which do not naturally conform to the standard input-state-output perspective from among physical systems, and notably the class of electric circuits, in addition to examples which contravene the typical controllability assumptions within the theory of dissipativity. 

In this article, we first review the celebrated theorems from the landmark Dissipative Dynamical Systems Parts I and II papers, in order to highlight their contributions and draw out outstanding issues and questions. These issues and questions are most naturally formulated and addressed using behavioral theory, and are brought into sharp focus by the consideration of the behaviors of electric circuits, as will be discussed in the subsequent two sections. We then present the recent results of the first author that address these questions, in the context of passive systems and electric circuit realization. Alongside the theory of dissipative dynamical systems, the landmark papers of Jan Willems also consider additional physical properties possessed by the classes of lossless, reciprocal, reversible, and relaxation systems. Once again, examples arise in physical systems, and electric circuit examples serve to illustrate the strengths and limitations of the theory presented in that paper. We further present recent results of the first author addressing the limitations of the theory of reciprocal systems, and discuss the consequences of these results in the context of lossless, reversible, and relaxation systems.

While the focus of this paper is expository, using electric circuit examples in part for their tutorial value, it is appropriate to note that the results presented herein are of equal applicability to the analysis and design of passive mechanical controllers, by drawing on the mechanical-electrical analogy. Indeed, the invention of the inerter \cite{mcs02} has motivated considerable recent research into passive network analysis and synthesis for mechanical applications such as the design of suspension systems for bridges, buildings, and automotive and railway vehicles (see the sidebar titled ``\nameref{sb:mn}'').


For the purposes of readability and accessibility, we will aim to use simple and minimal notation throughout. A considerable amount of the presented theory relies upon the algebra of matrices of real polynomials and rational functions in a single indeterminate, and the algebra of complex functions. Throughout, we will use the symbol $s$ to denote the indeterminate in a polynomial or rational function, and $\lambda$ to denote a complex number at which the polynomial is evaluated. For example, in the case of the scalar polynomial $a(s) = a_ns^n + a_{n-1}s^{n-1} + \ldots + a_1 s + a_0$ for real numbers $a_n, a_{n-1}, \ldots a_1, a_0$, and the complex number $\lambda$, we use the notation $a(s)$ to denote the polynomial itself, and $a(\lambda)$ to denote the complex number $a_n\lambda^n + a_{n-1}\lambda^{n-1} + \ldots + a_1 \lambda + a_0$. We say that a complex number $\lambda$ is in the closed right-half (respectively, closed left-half) plane if its real part is greater than or equal to zero (respectively, less than or equal to zero), and in the open right-half (respectively, left-half) plane if, in addition, its real part is non-zero. We further use a bar to denote the complex conjugate of a complex matrix, an asterisk to denote the Hermitian transpose of a complex matrix, and a superscript $T$ to denote its transpose. For example, for the matrix
$$
A = \begin{bmatrix}a_{11}& a_{12}\\ a_{21}& a_{22}\end{bmatrix},
$$
with $a_{11}, a_{12}, a_{21}$, and $a_{22}$ all complex numbers, then
$$
A^* = \bar{A}^T = \begin{bmatrix}\bar{a}_{11}& \bar{a}_{21}\\ \bar{a}_{12}& \bar{a}_{22}\end{bmatrix}.
$$
Furthermore, for a real and symmetric matrix $A$ (that is, $A = A^T$), then if $x^T A x \geq 0$ for all real vectors $x$ we say that $A$ is non-negative definite (denoted by $A \geq 0$), and if in addition $x^T A x = 0$ implies that $x$ is zero then we say that $A$ is positive definite (denoted $A > 0$). Additionally, for a pair of real and symmetric matrices $A$ and $B$, then $A \geq B$ denotes that $A-B \geq 0$.

\subsection{Dissipative Dynamical Systems Part 1}
In Dissipative Dynamical Systems Part 1 \cite{JWDSP1}, the concept of a dissipative system is defined in a very general setting, alongside associated concepts such as the available storage and the required supply, and fundamental results concerning stability and interconnection of dissipative systems. This unified previous work on optimal control and the passivity theorem and provided a natural generalisation of Lyapunov methods from autonomous to open dynamical systems. The concepts have stood the test of time, and inspired ongoing research into applications that include the analysis and design of interconnected (such as cyber-physical) systems; the theory of integral quadratic constraints and the related absolute stability problem; the concepts of stabilisation by adding dissipativity and passivity-based control; port Hamiltonian systems; network analysis and synthesis for electrical, mechanical and multiphysics systems; and $H_{\infty}$ methods. It is not possible to do justice to this large body of subsequent literature in this article, and we will only elaborate further on the final two applications and their relation to the results presented in this article (see the sidebars titled ``\nameref{sb:pna}'', ``\nameref{sb:mn}'', and ``\nameref{sb:hinf}''). In this section, we overview the key concepts and contributions of \cite{JWDSP1}, and draw out those assumptions and associated issues as recognised by Jan Willems himself (see, for example, \cite{JWDDS}), whose resolution we will discuss in the remainder of this article.

In \cite{JWDSP1}, the concept of a dissipative system is formalised as a system possessing a so-called \emph{storage function} with respect to a given \emph{supply rate}. This generalises the notion of a passive system such as an electric circuit, where the supply rate corresponds to the instantaneous power provided to the system. Two additional functions of the state of the system are defined---termed the \emph{available storage} and \emph{required supply}---which themselves turn out to be storage functions. These characterise the solutions to optimal control problems: those of maximising the future supply extracted from the system as a function of its initial state, and minimising the past supply input to the system in order to drive the system from a state of minimum storage to a given state. The latter concept of the required supply is well-posed only when there exists such a state of minimum storage, and only for those states that are reachable from such a minimum storage state.

The concept 
is developed in the context of a general stationary (that is, time-invariant) input-state-output system. Stationarity is assumed principally for notational convenience. The system is characterised in terms of an input $u$, state $x$, and output $y$ related via a state transition map $\phi$ and a read-out map $r$. Here, the state transition map is a function of the initial state $x_0$ and time $t_0$, the input applied $u$, and a time $t \geq t_0$, whereupon the state $x(t)$ at time $t$ satisfies 
\begin{equation}
\label{eq:xstf}
x(t) = \phi(t,t_0,x_0,u) \text{ for } t \geq t_0.
\end{equation}
The state transition map satisfies the natural axioms, as are satisfied by differential equations of the form $\dot{x} = f(x,u)$, of 1) consistency; 2) determinism; 3) the semi-group property; and 4)  stationarity. Here, axioms 1)--3) reflect the notion that the state at any given instant $t_1$ is uniquely determined by the state at any given prior instant $t_0 \leq t_1$ and the input in the interval from $t_0$ to $t_1$, and axiom 4) captures the time-invariance of the system. These four axioms are
characterised mathematically as the following four conditions,
which hold for any given initial state $x_0$, input $u$, and $t_2
\geq t_1 \geq t_0$:
\begin{enumerate}
\item $\phi(t_0, t_0, x_0, u) = x_0$; 
\item $\phi(t_1, t_0, x_0, u_1) = \phi(t_1, t_0, x_0, u_2)$ whenever $u_1(t) = u_2(t)$ for all $t_0 \leq t \leq t_1$; 
\item $\phi(t_2, t_0, x_0, u) = \phi(t_2, t_1, \phi(t_1, t_0, x_0, u), u)$; and 
\item if $\hat{u}(t) = u(t+T)$ for all real $t$, then $\phi(t_1+T,t_0+T,x_0,\hat{u})=\phi(t_1,t_0,x_0,u)$.
\end{enumerate}
The output map is then a function of the state $x$ and the applied input $u$, and determines the output $y$ via the relationship 
\begin{equation}
\label{eq:yrof}
y(t) = r(x(t),u(t)) \text{ for } t \geq t_0.
\end{equation}
Associated with the system is a supply rate, denoted $w$, which is taken to be a locally integrable function of the system's input and output. The dissipativity concept is then stated in Definition 2 of \cite{JWDSP1} in terms of the existence of a storage function. Specifically, a dynamical system is said to be dissipative with respect to a supply rate $w$ if there exists a non-negative function $S$ of the state of the system that satisfies
\begin{equation}
\label{eq:di}
S(x_0) + \int_{t_0}^{t_1}w(u(t),y(t))\mathrm{dt} \geq S(x(t_1))
\end{equation}
for all possible initial states $x_0 = x(t_0)$, times $t_1 \geq t_0$ and inputs $u$. In words, the system is dissipative if the change in storage from time $t_0$ to time $t_1$ (that is, $S(x(t_1)) - S(x_0)$) cannot exceed the total supply to the system in this interval.

The losslessness property is subsequently characterised in Definition 6 of the paper in terms of systems for which the change in storage is always equal to the total supply. Specifically, a dissipative dynamical system with supply rate $w$ and storage function $S$ is lossless if
$$
S(x_0) + \int_{t_0}^{t_1}w(u(t),y(t))\mathrm{dt} = S(x(t_1))
$$
for all possible initial states $x_0 = x(t_0)$ and inputs $u$.

Alongside the notions of dissipativity, 
the paper defines two functions of state: the \emph{available storage} (denoted $S_a$) and \emph{required supply} (denoted $S_r$). Here, for a given initial state $x_0$ at time $t_0$, the available storage $S_a(x_0)$ is the supremum of the extracted supply in the interval $t_0$ to $t_1$ (that is, $-\int_{t_0}^{t_1}w(u(t),y(t))\mathrm{dt}$) over all possible intervals $t_1 \geq t_0$ and inputs $u$. 

Of course, it may be the case that the supremum does not exist (that is, it is not finite). Indeed, Theorem 1 of \cite{JWDSP1} establishes that the available storage is finite for every given initial state $x_0$ if and only if the system is dissipative.

The required supply is defined on the assumption that there exists a state $x^*$ of minimum storage, and the storage function is normalised such that $S(x^*) = 0$. Then, for a given state $x_0$ at time $t_0$, the required supply $S_r(x_0)$ is the infimum of the provided supply in the interval $t_{-1}$ to $t_0$ (that is, $\int_{t_{-1}}^{t_0}w(u(t),y(t))\mathrm{dt}$) over all possible intervals $t_{-1} \leq t_0$ and inputs $u$ that drive the state from the state of minimum storage $x(t_{-1}) = x^*$ to $x(t_0) = x_0$.

Once again, it may be the case that the infimum does not exist (that is, it is not finite). Theorem 2 of  \cite{JWDSP1} establishes that the required supply is finite for every given state $x(t_0) = x_0$ if and only if the system is dissipative. Moreover, any given storage function $S$ is bounded between the available storage and required supply, that is, $S_a(x_0) \leq S(x_0) \leq S_r(x_0)$ for all initial states $x_0$. Significantly, \emph{these results are only valid for systems whose state space is reachable from the point of minimum storage}, that is, for any given point $x_0$ in the state space and any given time $t_0$, there exists a $t_{-1} \leq t_0$ and an input $u$ that drives the system from $x(t_{-1}) = x^*$ to $x(t_0) = x_0$. 

The theory proceeds through four further theorem statements, as follows:
\begin{enumerate}
\item Theorem 3 states that the set of possible storage functions forms a convex set. This follows in a straightforward manner from the definition of a dissipative system.
\item Theorem 4 states that if the system is reachable from and controllable to the state of minimum storage and is lossless, then the available storage and required supply are equal and hence the storage function is unique.
\item Theorem 5 shows that a neutral interconnection of dissipative systems is itself dissipative. A formal definition of interconnection and neutral interconnections is provided in Section 4 of \cite{JWDSP2} and relies on an assumption of well-posedness of the interconnection. It generalises the fact that passive systems interconnected in accordance with Kirchhoff's current and voltage laws are themselves passive. However, for passive systems such as electric circuits, it need not be the case that the interconnection is well-posed (see the sidebar titled ``\nameref{sb:pna}'').
\item Theorem 6 proves that the point of minimum storage ($x^*$) is a stable equilibrium point (but not necessarily asymptotically stable) if the storage function is continuous and attains a strong local minimum at $x^*$.
\end{enumerate}

While the conceptual framework introduced in \cite{JWDSP1} is very general, there remain certain limitations which, in particular, prevent the theory being universally applicable to certain physical systems such as electric circuits. Specifically, the underlying definition of a dissipative system relies on a representation of the system in terms of an input, state and output related via a state transition and read-out map. As emphasised in \cite{JWBAOIS, JWDDS}, many physical systems are not naturally characterised in such a form. More significantly, from an external perspective (that is, in terms of the system's inputs and outputs), the representation is non-unique, and the developed theory does not address the question of whether or not the dissipativity property depends on the specific representation. Indeed, in the context of so-called cyclo-dissipative systems, the dissipativity property has been shown to be representation dependent, particularly when one allows for unobservable states (see the sidebar titled ``\nameref{sb:cdac}'').

In the section titled ``\nameref{sec:psna}'', we present a trajectory-level (representation-free) definition of passivity that overcomes these limitations, is fitting for capturing the realizable behaviors of electric circuits, and can readily be extended to a more general definition of dissipativity. This is further developed within the context of finite dimensional linear time-invariant systems, thereby addressing some of the limitations of the theory presented in Part 2 of Jan Willems' Dissipative Dynamical Systems papers, which we discuss in the next section.


\subsection{Dissipative Dynamical Systems Part 2}

Part 2 of Dissipative Dynamical Systems \cite{JWDSP2} draws upon the network synthesis work of \cite{Kalman_newchar, YoulaTissi, AndVong}, the famous Kalman Yakubovich Popov lemma (or positive-real lemma), and linear optimal control, and unifies these under the dissipativity formalism. Specifically, \cite{JWDSP2} further develops the dissipativity concept for the case of finite dimensional linear time-invariant (FDLTI) systems where the relationship between input $u$, state $x$ and output $y$ takes the form
\begin{equation}
\dot{x} = Ax + Bu, \hspace{0.3cm} y = Cx + Du, \label{eq:ssr}
\end{equation}
for some real-valued matrices $A, B, C$ and $D$, and the supply rate takes the form
$$
w = u^T y.
$$
Furthermore, the system is assumed to be both controllable and observable. This represents the class of passive and controllable systems that can be represented in such an input-output form. The concept of the reciprocity of such a system is also introduced, and particular attention is paid to systems that are both passive and reciprocal. As discussed in the section titled ``\nameref{sec:pec}'', examples of such linear passive systems include electric circuits, which are also reciprocal when they contain only resistors, transformers, capacitors and inductors, and the developed theorems lead to network realizations of both passive systems and passive and reciprocal systems subject once again to the controllability caveat and the input-output formalism. The generalisation to arbitrary quadratic supply rates is then discussed in \cite[Remark 9]{JWDSP2}. 

Theorem 1 of \cite{JWDSP2} presents the classical frequency-domain condition for dissipativity in the case of the supply rate $w = u^T y$. Namely, for a minimal (that is, controllable and observable) input-state-output system taking the form of equation \eqref{eq:ssr}, the system is dissipative with respect to the supply rate $w = u^T y$ if and only if the transfer function $G(s) = D + C(sI - A)^{-1}B$ is \emph{positive-real}, that is, $G(\lambda) + G(\lambda)^* \geq 0$ for all $\lambda$ in the open right-half plane.

In Lemma 2 and Theorem 2 of \cite{JWDSP2}, dissipativity is characterised for those minimal input-state-output systems taking the form of equation \eqref{eq:ssr} with $D+D^T$ invertible in terms of the existence of solutions to the algebraic Riccati equation
\begin{equation}
XA + A^T X + (XB-C^T)(D+D^T)^{-1}(B^TX - C) = 0.\label{eq:are}
\end{equation}
Specifically, the system is dissipative with respect to the supply rate $w = u^T y$ if and only if there exists an $X \geq 0$ satisfying this algebraic Riccati equation. Furthermore, the available storage satisfies $S_a(x) = \tfrac{1}{2}x^TX^{-}x$ where $X^{-}$ is the unique real symmetric matrix which satisfies the algebraic Riccati equation and for which the eigenvalues of $A+B(D+D^T)^{-1}(B^TX^{-} - C)$ are all in the closed left-half plane. In addition, the required supply satisfies $S_r(x) = \tfrac{1}{2}x^TX^{+}x$ where $X^{+}$ is the unique real symmetric matrix which satisfies the algebraic Riccati equation and for which the eigenvalues of $A+B(D+D^T)^{-1}(B^TX^{+} - C)$ are all in the closed right-half plane. 

Theorem 3 of \cite{JWDSP2} then characterises dissipativity for cases in which $D+D^T$ need not be invertible in terms of the matrix inequality
\begin{equation}
\begin{bmatrix} -A^T X - XA& C^T - XB\\ C - B^T X& D + D^T\end{bmatrix} \geq 0. \label{eq:lmi}
\end{equation}
Specifically, for a minimal input-state-output system taking the form of \eqref{eq:ssr}, the system is dissipative if and only if there exists an $X \geq 0$ satisfying the above matrix inequality. In particular, the available storage and required supply take the form $S_a(x) = \tfrac{1}{2}x^T X^- x$ and $S_r(x) = \tfrac{1}{2}x^T X^+ x$ where $X^- > 0$ and $X^+ > 0$ are solutions to the matrix inequality, and any other non-negative definite solution to the matrix inequality satisfies $X^+ \geq X \geq X^-$.

Then, Theorem 4 of \cite{JWDSP2} characterises the dissipation of a dissipative system in terms of factorisations of the matrix in inequality \eqref{eq:lmi} of the form
\begin{equation}
\begin{bmatrix} -A^T X - XA& C^T - XB\\ C - B^T X& D + D^T\end{bmatrix} = \begin{bmatrix} M^T\\ N^T\end{bmatrix}\begin{bmatrix} M& N\end{bmatrix}. \label{eq:dr}
\end{equation}
Such a factorisation exists if and only if the matrix inequality \eqref{eq:lmi} holds. Moreover, for any such factorisation, the Euclidean norm of $Mx + Nu$ characterises the dissipation rate of the system, in the sense that the change in storage between states $x(t_0)$ at time $t_0$ and $x(t_1)$ at time $t_1$ satisfies
\begin{multline*}
\frac{1}{2}x(t_0)^TXx(t_0) - \frac{1}{2}x(t_1)^TXx(t_1) \\= -\int_{t_0}^{t_1}u(t)^Ty(t)\mathrm{dt} + \frac{1}{2}\int_{t_0}^{t_1}(Mx(t) + Nu(t))^T(Mx(t) + Nu(t))\mathrm{dt}.
\end{multline*}

Theorem 5 of \cite{JWDSP2} then provides a frequency domain condition for a system to be lossless with respect to the supply rate $u^Ty$, namely that the transfer function $G(s) = D+C(sI-A)^{-1}B$ be lossless positive-real, that is, $G(s)$ is positive-real and $G(i\omega) + G(-i\omega)^T = 0$ for all real $\omega$ other than imaginary axis poles of $G(s)$. It is further shown that the system is lossless if and only if the matrix inequality \eqref{eq:lmi} has a unique solution, and for this unique solution the inequality holds with equality.

The paper \cite{JWDSP2} goes on to characterise reciprocal systems in terms of those with a signature symmetric transfer function. Specifically, a signature matrix refers to a diagonal matrix whose diagonal entries are all either $+1$ or $-1$, and a system taking the form of \eqref{eq:ssr} is called reciprocal with respect to the external signature matrix $\Sigma_e$ if $\Sigma_e(D + C(sI - A)^{-1}B)$ is symmetric. 

In Theorem 6 of \cite{JWDSP2}, it is shown that if a minimal input-state-output system taking the form of equation \eqref{eq:ssr} is reciprocal with respect to the external signature matrix $\Sigma_e$, then the system has a minimal realization
\begin{align}
\dot{z} &= A_1 z + B_1 u, \label{eq:assr} \\
y &= C_1 z + D_1 u, \label{eq:assr2}
\end{align}
where
\begin{equation}
\begin{bmatrix}\Sigma_i& 0\\ 0& \Sigma_e\end{bmatrix}\begin{bmatrix}-A_1& -B_1\\ C_1& D_1\end{bmatrix} \text{ is symmetric} \label{eq:ir}
\end{equation}
for some signature matrix $\Sigma_i$. Moreover, any system admitting such a realization is reciprocal with respect to the signature matrix $\Sigma_e$. We call a realization of the aforementioned form an \emph{internally reciprocal} realization, and we refer to those states corresponding to a $+1$ in the internal signature matrix $\Sigma_i$ as states of \emph{even parity} and those corresponding to a $-1$ as states of \emph{odd parity}.

In Theorem 7 of \cite{JWDSP2}, it is shown that if a minimal input-state-output system taking the form of equation \eqref{eq:ssr} is dissipative with respect to the supply rate $w = u^T y$ and reciprocal with respect to the external signature matrix $\Sigma_e$, then the system admits a minimal realization taking the form of \eqref{eq:assr}--\eqref{eq:assr2}, where \eqref{eq:ir} holds for some signature matrix $\Sigma_i$, and also
\begin{equation}
\begin{bmatrix}-A_1& -B_1\\ C_1& D_1\end{bmatrix} + \begin{bmatrix}-A_1& -B_1\\ C_1& D_1\end{bmatrix}^T \geq 0. \label{eq:ipir1}
\end{equation}
We refer to such a realization as an \emph{internally reciprocal and passive} realization. The terminology here is motivated in part by the fact that such a realization leads to a realization of the controllable input-output behavior as the driving-point behavior of an electric circuit comprising the passive and reciprocal elements: resistors, inductors, capacitors, and transformers (see the sidebar titled ``\nameref{sb:pns}'').

The paper concludes with two final theorems on reversible systems and relaxation systems. Here, a system is said to be \emph{reversible} with respect to a signature matrix $\Sigma_e$ if, whenever two inputs $u_1$ and $u_2$ satisfy $u_2(t) = -\Sigma_e u_1(-t)$, and the corresponding states of the system are zero when $t = 0$, then the corresponding outputs $y_1$ and $y_2$ satisfy $y_2(t) = \Sigma_e y_1(-t)$. In Theorem 8 of \cite{JWDSP2} it is shown that, for a minimal input-state-output system that is dissipative with respect to the supply rate $u^T y$, reversibility with respect to a signature matrix $\Sigma_e$ is equivalent to being lossless and reciprocal with respect to that signature matrix. Thus, dissipative reversible systems correspond to electric circuits comprising inductors, capacitors, and transformers. 

A relaxation system is then defined as a system whose impulse response is completely monotonic and symmetric, that is, with $G(t) = Ce^{At}B$, then $D = D^T \geq 0$, $G(t) = G(t)^T$, $G(t) \geq 0$ for all $t \geq 0$, and $(-1)^n\tfrac{d^n}{dt^n}G(t) \geq 0$ for all integer $n$ and all $t \geq 0$. These arise as the driving-point voltage response of electric circuits comprising resistors, capacitors, and transformers to an impulse in the driving-point current, and the driving-point current response of circuits comprising resistors, inductors, and transformers to an impulse in the driving-point voltage. The final theorem of \cite{JWDSP2} shows that a minimal input-state-output system is a relaxation system if and only if it is passive and reciprocal with all states having even parity. In other words, the system has a minimal realization taking the form of equations \eqref{eq:assr}--\eqref{eq:assr2} that satisfies \eqref{eq:ipir1} and where \eqref{eq:ir} holds with $\Sigma_i$ being the identity matrix.

The principal limitations of the presented results in \cite{JWDSP2} are the input-state-output formalism and the assumptions of controllability and observability. Indeed, in the absence of controllability, then the positive-real condition is a necessary but not sufficient condition for passivity. This follows as positive-realness is a condition on the transfer function of a system, which fails to differentiate between, for example, the behaviors of the systems described by the differential equations (i) $y = u$, (ii) $\tfrac{dy}{dt} + y = \tfrac{du}{dt} + u$, and (iii) $\tfrac{dy}{dt} - y = \tfrac{du}{dt} - u$, all of whom have the transfer function $G(s) = 1$. Indeed, it can be shown that the first two such systems are passive but the third is not. For example, it admits the trajectory $y(t) = e^{t}$ and $u(t) = -e^{t}$, for which $-\int_{t_0}^{t_1}u(t)y(t)\mathrm{dt} = \tfrac{1}{2}(e^{2t_1} - e^{2t_0})$, and this can be made arbitrarily large by choosing a sufficiently large $t_1$. Moreover, the very definitions of reciprocal, reversible, and relaxation systems provided in \cite{JWDSP2} are formulated in terms of input-output systems, yet these concepts serve to describe the behaviors of important subclasses of electric circuits that need not naturally conform to such an input-output framework. 

The aforementioned issues are clarified by introducing a trajectory-level and representation-free definition of passivity and also of reciprocity, and by characterising systems in terms of the differential equation governing their behavior instead of focussing on input-output properties such as the transfer function or impulse response. Accordingly, we will proceed to outline a theory of passivity and reciprocity for systems corresponding to the set of solutions to linear differential equations (as opposed to systems formulated using an input-state-output representation), and the consequences of such results in the context of reversible and relaxation systems. Specifically, we will report on the following:
\begin{enumerate}
\item An extension of the positive-real condition to provide a necessary \emph{and sufficient} condition for passivity for systems that need not be controllable.
\item Characterisations for passivity in terms of the solutions to the algebraic Riccati equation \eqref{eq:are} and the matrix inequality \eqref{eq:dr} that hold in the absence of controllability or observability assumptions.
\item Necessary and sufficient conditions for the realizability of the solutions to a given differential equation as the driving-point behavior of a passive electric circuit.
\item An extension of the results on signature-symmetric realizations to systems that need not be controllable.
\item Necessary and sufficient conditions for the realizability of a behavior as the driving-point behavior of a passive and reciprocal electric circuit.
\item The implications of these results for lossless systems, reversible systems, and relaxation systems.
\end{enumerate} 

\section{Passive electric circuits}
\label{sec:pec}
A classical example of a dissipative system is a passive electric circuit (see the sidebar titled ``\nameref{sb:pna}''). Yet such circuits do not immediately lend themselves to an input-output or input-state-output description. Instead, they are characterised in terms of pairs of driving-point variables---voltage and current---whose inner product determines the instantaneous power provided to the circuit. In many cases a given electric circuit may be considered as an input-state-output system with the driving-point currents as input, the inductor currents and capacitor voltages as states, and the driving-point voltages as output. However, this is not always possible, and such an input-state-output description is not a natural starting point for the modelling of such circuits. Indeed, it is not necessarily the case that the vector of driving-point currents or voltages can be considered as an input to the system. A simple example is the transformer, for which it is not possible to arbitrarily assign all of the port voltages, nor is it possible to arbitrarily assign all of the port currents.

Taking inspiration from the behavioral approach (see, for example, \cite{JWBAOIS}), a more natural starting point for the modelling of such a system is to write down the equations governing the system's full behavior itself---that is, the differential equations corresponding to the inductor, capacitor, resistor, transformer, and gyrator elements, and the algebraic interconnection equations corresponding to Kirchhoff's current and voltage laws. The internal element voltages and currents can then be eliminated in order to find a differential equation relating the driving-point currents and voltages. This can be handled systematically and formally using the element extraction approach, alongside the elimination method formalised in Jan Willems' behavioral theory, as outlined in \cite{pec_ba}. As noted by Jan Willems himself in \cite[Section 12]{JWDDS}, two questions then arise naturally: 
\begin{enumerate}
\item Can the set of differential equations realizable by passive electric circuits be characterised?
\item Given a `realizable' differential equation, is there a constructive procedure that obtains a passive electric circuit realization?
\end{enumerate}
Alongside these questions are two parallel questions relating to reciprocity:
\begin{enumerate}
\item Can the set of differential equations realizable by passive \emph{and reciprocal} electric circuits be characterised?
\item Given a `realizable' differential equation, is there a constructive procedure that obtains a passive \emph{and reciprocal} electric circuit realization?
\end{enumerate}

The answers to such questions of realizability in the case of controllable systems that possess a so-called \emph{hybrid matrix description} is classical and follows from the internally passive and reciprocal realization of \cite[Theorem 7]{JWDSP2}. A detailed account is provided in \cite{AndVong}. The general case (to include uncontrollable systems and those that are not naturally characterised in an input-output setting) was only recently solved (see \cite{THTPLSNA, THRSC, pec_ba}), and follows from a more general theory of passive systems that we describe in the section titled ``\nameref{sec:psna}''. This general theory is most naturally formalised using the framework of the behavioral approach of Jan Willems, and we summarise the relevant concepts and results from the behavioral approach in the next section.

\section{The behavioral approach}
\label{sec:ba}
The philosophy underlying the behavioral approach is that many systems do not immediately lend themselves to an input-output or input-state-output viewpoint, and to take such a viewpoint as the starting point for the modelling of general physical systems is conceptually limiting and potentially restrictive. Instead, the behavioral approach characterises systems as sets of trajectories, being vector-valued functions of time, taken from a function space termed the universum (such as the space of vector-valued locally integrable functions), and restricted to those elements of the universum that are consistent with the laws governing the system. Such laws are typically specified by equations, such as the constituent element laws and Kirchhoff's current and voltage laws governing the behavior of a given electric circuit. In the case of finite dimensional linear time-invariant (FDLTI) equations, the natural characterisation of such systems is the set of solutions to a differential equation of the general form
\begin{equation}
R_n \tfrac{d^n w}{dt^n} + R_{n-1}\tfrac{d^{n-1} w}{dt^{n-1}} + \ldots + R_1 \tfrac{dw}{dt} + R_0w = 0
\end{equation}
for some real-valued matrices $R_0, R_1, \ldots , R_n$, where the behavior of the system corresponds to the set of solutions $w$ to the differential equation, taken from a suitable function space. To allow for the consideration of, for example, step functions, while acknowledging that the dissipativity concept requires the evaluation of integrals of the products of the system's variables, then it is fitting to take the set of vector-valued locally square integrable functions as the function space. The differential equation is then to be treated in a weak sense in the manner formalised in \cite{JWIMTSC}. 

It is notationally and conceptually convenient to define the real polynomial matrix
\begin{equation*}
R(s) = R_n s^n + R_{n-1} s^{n-1} + \ldots + R_1 s + R_0,
\end{equation*}
whereupon the behavior can be considered as the kernel of the polynomial differential operator $R(\tfrac{d}{dt})$. It is then possible to analyse such linear systems using matrix polynomial algebra techniques. For example, such techniques allow for the elimination of variables as described in the sidebar titled ``\nameref{sb:etape}''. This is applicable to, for example, obtaining the driving-point behavior of an electric circuit by eliminating the internal currents and voltages from the set of equations obtained from Kirchhoff's laws and the constituent element laws.

While the behavioral approach emphasises the benefits of placing all variables on an equal footing by operating with high order differential equations as opposed to taking an input-state-output model from the outset, it is \emph{not} the case that the approach seeks to dispose with inputs, states, and outputs outright. Rather, the behavioral approach articulates the qualities of variables within a system that bestow them with the properties of inputs, states, and outputs. Here, the state of the system is characterised in terms of concatenability of trajectories---namely, two trajectories $w_1$ and $w_2$ are said to have the same \emph{state} at an instant $t_0$ if they can be concatenated at that instant, that is, if the trajectory $w$ is also in the behavior where $w$ satisfies $w(t) = w_1(t)$ for all $t < t_0$ and $w(t) = w_2(t)$ for all $t \geq t_0$ (see, for example, \cite[Section 2]{JWDSP1}, \cite[Section 3]{PRSMLS}). Then, a set of variables within a system comprise an \emph{input} $u$, with the remaining variables constituting an \emph{output} $y$, if the input variables are free (that is, they are unconstrained by the laws of the system) and the output variables do not contain any further free components \cite[Def. 3.3.1]{JWIMTSC}. For such an input-output partition, the state of the system at an instant $t_0$, together with the input $u(t)$ for $t \geq t_0$, determine the output $y(t)$ for all $t \geq t_0$. Indeed, for FDLTI systems, it is proven that it is always possible to partition the system into inputs and outputs, whereupon an input-state-output representation for the system can readily be obtained. In fact, it is always possible to construct an \emph{observable} input-state-output representation taking the form of \eqref{eq:ssr} (in the sense that the state of the system at a given time $t_0$ is uniquely determined by the system's input and output for $t \geq t_0$), in which case two trajectories $(u_1, y_1)$ and $(u_2, y_2)$ have the same state at time $t_0$ if and only if the corresponding trajectories $(u_1, y_1, x_1)$ and $(u_2, y_2, x_2)$ satisfy $x_1(t_0) = x_2(t_0)$. For example, the papers \cite{PRSMLS, AVDSSMIP, THBRSF} outline constructive techniques for obtaining observable input-state-output representations for any given behavior.

Another strength of the behavioral theory is the formalism of the concepts of controllability and stabilisability as trajectory-level properties as opposed to properties of a given input-state-output representation. Conceptually, a behavior is called controllable if for any given trajectories $w_1$ and $w_2$ from the behavior and any given real number $t_0$, there exists a $t_1 \geq t_0$ and a trajectory $w$ from the behavior where $w(t) = w_1(t)$ for all $t < t_0$ and $w(t) = w_2(t)$ for all $t \geq t_1$. Similarly, a behavior is called stabilisable if for any given trajectory $w_1$ from the behavior and any given real number $t_0$, there exists a trajectory $w$ from the behavior with $w(t) = w_1(t)$ for all $t < t_1$ and $\lim_{t \rightarrow \infty} w(t) = 0$. It is appropriate to note that any input-state-output representation of a given uncontrollable (resp.\ unstabilisable) input-output behavior is itself uncontrollable (resp.\ unstabilisable). The converse also holds whenever the input-state-output representation is observable, but does not hold in general (see \cite[Note D4]{THTPLSNA}). 

Algebraically, the behavior corresponding to the kernel of the polynomial differential operator $R(\tfrac{d}{dt})$ is controllable if and only if the rank of $R(\lambda)$ has the same value for all complex numbers $\lambda$, and stabilisable if and only if the rank of $R(\lambda)$ has the same value for all $\lambda$ in the closed right-half plane \cite{JWIMTSC}. This highlights a distinct advantage of the polynomial differential operator representation relative to characterising systems in terms of their transfer function. For example, the single-input single-output system $\tfrac{dy}{dt} + y = \tfrac{du}{dt} + u$ permits all trajectories of the form $y(t) = u(t) + ke^{-t}$ with $k$ a real number. This behavior corresponds to the kernel of the differential operator $R(\tfrac{d}{dt})$ with $R(s) = \begin{bmatrix}s+1& -(s+1)\end{bmatrix}$. Since $R(-1)= \begin{bmatrix}0&0\end{bmatrix}$, and the rank of $R(\lambda)$ is one for all $\lambda\neq -1$, this behavior is not (behaviorally) controllable. Indeed, the transfer function of this system from $u$ to $y$ is equal to $1$ (owing to a cancelling factor of $s+1$), and only identifies the trajectories $y(t) = u(t)$. In the context of passivity, the transfer function condition of positive-realness does not capture the constraints on such cancelling factors that may arise as a result of the system's passivity. This naturally leads to the following question, first posed in \cite[Section 12]{JWDDS}: what are the necessary and sufficient conditions on these factors that guarantee the passivity of the system? The answer to this question was provided in \cite{THTPLSNA} and is summarised in the next section of this article.

At the opposite extreme to \emph{controllable} behaviors are \emph{autonomous} behaviors, which generalise the notion of the autonomous state-space system $\tfrac{dx}{dt} =Ax$ to systems described by high order differential equations. Specifically, an \emph{autonomous} behavior is one for which the future trajectory (for $t \geq t_0$ for any given real number $t_0$) is fully determined by the state of the system at $t_0$. Algebraically, such systems correspond to the solutions to differential equations of the form $R(\tfrac{d}{dt})w = 0$ where $R(s)$ is a square matrix whose determinant $\det{(R(s))}$ is non-zero. In addition, the behavior of such an autonomous system is a finite dimensional linear subspace of the set of vector-valued infinitely differentiable functions, whose dimension is equal to the degree of $\det{(R(s))}$, the trajectories being linear sums of vectors of products of polynomials, sinusoids, and exponential functions of time (see \cite[Theorem 3.2.16]{JWIMTSC}). 

It can further be shown that the behavior of any given FDLTI system can be decomposed into a \emph{controllable} part and an \emph{autonomous} part. Specifically, there exist controllable and autonomous sub-behaviors such that any trajectory $w$ from the originating behavior can be uniquely decomposed as a sum $w = w_1 + w_2$ where $w_1$ belongs to the controllable sub-behavior and $w_2$ belongs to the autonomous sub-behavior. Such sub-behaviors can be obtained using polynomial matrix algebra techniques, see for example, \cite[Lemma 11]{THRSC}. We then refer to the number of uncontrollable modes of the system as the dimension of the autonomous sub-behavior, and we note that this is equal to the number of uncontrollable modes in any observable input-state-output representation for the behavior.

\section{Passive systems}
\label{sec:psna}
The aforementioned discussion of electric circuits and the theory of behaviors motivates a theory of passive systems that
\begin{enumerate}
\item Does not take an input-state-output description as its starting point, but is instead characterised by pairs of variables (such as the driving-point voltages and currents of an electric circuit) whose inner product is the supply rate (the instantaneous power supplied to the system), and which in the case of FDLTI systems are characterised by the set of solutions to a linear differential equation.
\item Does not assume the system to be controllable, or even to presume a-priori that the system is stabilisable.
\end{enumerate}

The starting point for the development is a differential equation of the form
\begin{equation}
P(\tfrac{d}{dt})i = Q(\tfrac{d}{dt})v, \label{eq:gcfps}
\end{equation}
where $P(s)$ and $Q(s)$ are square polynomial matrices, and the rows of the matrix $\begin{bmatrix}P(\lambda)& {-}Q(\lambda)\end{bmatrix}$ are independent for almost all complex values of $\lambda$. The driving-point behavior of any given electric circuit always takes such a form \cite{pec_ba}. There is no a-priori assumption here that either $i$ or $v$ have the property of an input, and indeed it may be the case that neither $i$ nor $v$ have this property (we recall the prior example of a transformer). 

\subsection{Passivity}
\label{sec:passivity}
We say that a system taking the form of equation \eqref{eq:gcfps} is passive if, for any given real number $t_0$ and any given state of the system at $t_0$, the future extracted energy
\begin{equation}
-\int_{t_0}^{t_1}i(t)^Tv(t)\mathrm{dt}
\end{equation}
is bounded above for all $t_1 \geq t_0$ and for all future trajectories of the system that share the specified state at $t_0$. Keeping in mind the existence of an observable input-state-output representation taking the form of \eqref{eq:ssr} for any given FDLTI system, where $(u,y)$ corresponds to some partitioning of the variables $(i,v)$, then this is equivalent to the boundedness of the available storage $S_a(x_0)$ as defined in \cite{JWDSP1} for all real vectors $x_0$, where the supply rate $w$ is a quadratic function such that $i^Tv = w(u,y)$. In fact, for a passive system, there always exists an input-output transition such that $w(u,y) = u^Ty$ \cite[Theorem 9]{THTPLSNA}. Note, however, that this need not be the case for systems that are not passive. For instance, consider the system 

\begin{equation*}
\begin{bmatrix}0 & 1 \\0 & 0\end{bmatrix}\begin{bmatrix}i_1\\i_2\end{bmatrix} = \begin{bmatrix}0 & 0 \\0 & 1\end{bmatrix}\begin{bmatrix}v_1\\v_2\end{bmatrix},
\end{equation*}  
which is in the form of \eqref{eq:gcfps}. 
For this system, $i_2 = v_2 = 0$ so the only free variables that can form the input are $i_1$ and $v_1$. Thus, if $(u,y)$ were an input output partition then $u^Ty$ must take the form $i_1i_2 + v_1v_2$ or $i_1v_2 + v_1i_2$, which in neither case is equal to $i^Tv$.

We note that this trajectory-level definition of passivity does not depend in any crucial way on the linearity, time-invariance, or finite dimensionality of the system and it is thus fitting as a general definition suitable for non-linear, time varying, or infinite dimensional systems. It can also be naturally extended to a general trajectory-level definition of dissipativity. Indeed, in \cite{Hughes2018}, this trajectory-level definition was used to study so-called non-expansive systems, corresponding to input-state-output systems, as in equation \eqref{eq:ssr}, that are dissipative with respect to the supply rate $w(u,y) = u^T u - y^T y$. This led to results analogous to those detailed in this article that extend the classical bounded-real lemma to input-state-output systems that need not be controllable or observable (see the sidebar titled ``\nameref{sb:hinf}''). 

\subsection{Reciprocity}
\label{sec:recip}
Alongside the aforementioned trajectory-level definition of passivity, it is appropriate to consider a trajectory-level definition of reciprocity, and here we turn to \cite[Definition 2.7]{newclms}. Specifically, we say that a system taking the form of equation \eqref{eq:gcfps} is \emph{reciprocal} if, for any two trajectories $(i_a, v_a)$ and $(i_b, v_b)$ from the behavior that start at rest (that is, for which there exists a real number $t_0$ such that $i_a(t) = i_b(t) = v_a(t) = v_b(t) = 0$ for all $t < t_0$), then $\int_{-\infty}^{\infty}v_b(\tau)^Ti_a(t-\tau)\mathrm{d\tau} = \int_{-\infty}^{\infty}i_b(\tau)^Tv_a(t-\tau)\mathrm{d\tau}$ for all real $t$. We note that this is a property of the controllable part of the behavior alone. We do not intend to enter here into discussions about the physical origins of reciprocity in order to justify such a definition. Instead, we take as sufficient justification the fact that, in the FDLTI case, it corresponds precisely to the achievable external behaviors of systems possessing signature-symmetric realizations, and also (in the passive case) to the realizable driving-point behaviors of electric circuits comprising resistors, inductors, capacitors, and transformers.

We note that if the variables $i$ have the property of an input to the system, then the transfer function from $i$ to $v$ is symmetric. This is consistent with the system being reciprocal with respect to the identity matrix, in the sense of \cite{JWDSP2}. Similarly, if there exists a subset of the variables $i$ together with a complementary subset of the variables $v$ that constitute an input to the system, then the resulting input-output system is reciprocal with respect to the signature matrix $\Sigma_e$, in the sense defined in \cite{JWDSP2}, where $\Sigma_e$ is a signature matrix whose $k$th diagonal entry is $+1$ when the $k$th input corresponds to the $k$th variable from $i$, and $-1$ when it corresponds to the $k$th variable from $v$. Note that this definition implies that not all systems with signature-symmetric transfer functions are reciprocal. An example here being the gyrator, whose driving-point behavior is shown in Table \ref{fig:ece} in the sidebar titled ``\nameref{sb:pna}''. This is because an interconnection of passive elements including gyrators will not generally possess a signature-symmetric transfer function. In contrast, it can be shown that an interconnection of elements that are reciprocal in the sense defined in this subsection will itself be reciprocal whenever the interconnection laws correspond to Kirchhoff's current and voltage laws.

\subsection{Losslessness}
We call a behavior lossless if, for any given real numbers $t_1 \geq t_0$, any given states $x_0$ and $x_1$, and any given trajectory $(i, v)$ from the behavior that carries the system from state $x_0$ at time $t_0$ to state $x_1$ at time $t_1$, the available storage at state $x_0$ plus the supply to the system in the interval from $t_0$ to $t_1$ (that is, $\int_{t_0}^{t_1}i(t)^Tv(t)\mathrm{dt}$) is equal to the available storage at state $x_1$. Thus, the energy extracted from the system in the interval from $t_0$ to $t_1$ depends only on the change in the state from $x_0$ at time $t_0$ to $x_1$ at time $t_1$ and not on the specific trajectory that carries the system between these two states. Here, we again adopt the representation-free definition of state whereupon this definition is itself representation-free, where $x_0$ and $x_1$ denote equivalence classes amongst sets of trajectories that are concatenable at the instants $t_0$ and $t_1$ respectively, and we recall that these can be associated with real vectors corresponding to the values of $x(t_0)$ and $x(t_1)$ in any observable input-state-output representation for the behavior. It is also pertinent for the subsequent discussion to introduce the concept of a lossless sub-behavior of a given passive behavior. Here, we say a behavior $\mathcal{B}_1$ is a lossless sub-behavior of a given passive behavior $\mathcal{B}$ if $\mathcal{B}_1$ is controllable, any trajectory in $\mathcal{B}_1$ is also in $\mathcal{B}$, and the sub-behavior $\mathcal{B}_1$ is lossless.

\subsection{Necessary and sufficient conditions for passivity}
In the case of a controllable behavior in input-output form, the behavior is passive if and only if its transfer function is \emph{positive-real}. In the case of uncontrollable behaviors, it is also necessary to specify the constraints that passivity imposes on the autonomous part of the behavior. Such constraints were established in \cite{THTPLSNA}. Indeed, it is shown prior to Lemma 18 on p.\ 92 of \cite{THTPLSNA} that, for not necessarily controllable behaviors in the form of equation \eqref{eq:gcfps}, the behavior $\mathcal{B}$ is passive if and only if the following three conditions hold:
\begin{enumerate}
\item The controllable part of the behavior is passive.
\item $\mathcal{B}$ is (behaviorally) stabilisable (equivalently, the autonomous part of the behavior is stable).
\item If $(i_1,v_1)$ is a trajectory from a lossless sub-behavior of $\mathcal{B}$ that is at rest at times $t_0$ and $t_1$ (that is, $i_1$ and $v_1$ are both identically zero in a neighbourhood of both $t_0$ and $t_1$), and $(i_2, v_2)$ is a trajectory from the autonomous part of $\mathcal{B}$, then $\int_{t_0}^{t_1}i_1(t)^T v_2(t) + i_2(t)^T v_1(t)\mathrm{dt} = 0$.
\end{enumerate}
The necessity of the first condition is immediate. 
The necessity of the second condition can be ascertained by showing that if the behavior is not behaviorally stabilisable then, for any given real numbers $K$ and $t_0$, there exists a trajectory $(i,v)$ and a $t_1 \geq t_0$ such that the energy extracted from the system in the interval $t_0$ to $t_1$ (that is, $-\int_{t_0}^{t_1}i(t)^Tv(t)\mathrm{dt}$) exceeds $K$. The necessity of the third condition follows by noting that $(i, v)$ is a trajectory where $i = \alpha i_1 + i_2$ and $v = \alpha v_1 + v_2$ for any given real number $\alpha$, and since $i_1$ and $v_1$ are both at rest at both $t_0$ and $t_1$, then the state of the trajectory is independent of the choice of $\alpha$. Then $-\int_{t_0}^{t_1}i(t)^Tv(t)\mathrm{dt} = -\alpha^2 \int_{t_0}^{t_1}i_1(t)^Tv_1(t)\mathrm{dt} - \alpha \int_{t_0}^{t_1}i_1(t)^T v_2(t) + i_2(t)^T v_1(t)\mathrm{dt} -\int_{t_0}^{t_1}i_2(t)^Tv_2(t)\mathrm{dt}$, where the first term is zero since $(i_1, v_1)$ is from a lossless sub-behavior and is at rest at both $t_0$ and $t_1$, and the second term can be made arbitrarily large by choice of $\alpha$. 

That the three aforementioned conditions are also sufficient to guarantee passivity is established in \cite{THTPLSNA}. Indeed, they lead to algebraic tests for passivity that extend the positive-real condition to behaviors that need not be controllable. Specifically, it is shown in the so-called \emph{passive behavior theorem} \cite[Theorem 9]{THTPLSNA} that the behavior in equation \eqref{eq:gcfps} is passive if and only if the following three conditions all hold:
\begin{enumerate}
\item $P(\lambda)Q(\bar{\lambda})^T + Q(\lambda)P(\bar{\lambda})^T \geq 0$ for all $\lambda$ in the closed right-half plane.
\item The rows of $\begin{bmatrix}P(\lambda)& {-}Q(\lambda)\end{bmatrix}$ are independent for all $\lambda$ from the closed right-half plane.
\item If $p$ is a polynomial vector and $\lambda$ a complex number that satisfy $p(s)^T(P(s)Q(-s)^T+Q(s)P(-s)^{T}) = 0$ and $p(\lambda)^T\begin{bmatrix}P(\lambda)& {-}Q(\lambda)\end{bmatrix} = 0$, then $p(\lambda) = 0$.
\end{enumerate}

The second and third conditions can be evaluated using techniques from polynomial matrix algebra that can be easily implemented using standard symbolic algebra software. Also, several equivalent conditions to the third are presented in \cite{THTPLSNA}. In addition, while it may be the case that neither $i$ nor $v$ have the property of an input, it is shown that there will always exist an input comprised of a subset of the variables in $i$ together with the complementary variables in $v$. Thus, for a passive system, there always exists an input $u$ with the property that the instantaneous power delivered to the system takes the form $u^Ty$. 

Taking into account the preceding discussion, it follows that any passive system can indeed always be represented by an input-state-output system, as in equation \eqref{eq:ssr}, that is dissipative with respect to the supply rate $u^T y$. In this context, it can then be shown that dissipativity of the system is equivalent to the existence of a real symmetric non-negative definite matrix $X$ that satisfies inequality \eqref{eq:lmi}, and moreover that the available storage satisfies $S_a(x_0) = x_0^T X x_0$ where $X$ is the least non-negative definite matrix that satisfies this inequality. This extends \cite[Theorem 3]{JWDSP2} to the case of input-state-output systems that are not necessarily controllable or observable. Furthermore, it is shown in \cite[Theorem 13]{Hughes2018} that, in terms of dissipation rates, this $X$ has the property that the matrices $M$ and $N$ in equation \eqref{eq:dr} are such that $Z(s) = N + M(sI-A)^{-1}B$ is a spectral factor for $G(s) + G(-s)^T$, that is, that $G(s) + G(-s)^T = Z(-s)^TZ(s)$, that $Z(s)$ is analytic in the open right-half plane, and that the rows of $Z(\lambda)$ are independent for all $\lambda$ in the open right-half plane. Note, in contrast to the controllable case, the set of solutions to the matrix inequality will not be bounded above for an uncontrollable system owing to the ill-posedness of the required supply in this case (see \cite[Remark 14]{Hughes2018}).

In the case that $D+D^T$ is nonsingular, then it is also shown in \cite[Theorem 12]{Hughes2018} that the passivity of a system is equivalent to the existence of a non-negative definite matrix $X$ such that the algebraic Riccati equation \eqref{eq:are} holds. Moreover, the available storage satisfies $S_a(x_0) = x_0^T X x_0$ where $X$ is a solution to this algebraic Riccati equation for which the eigenvalues of $A+B(D+D^T)^{-1}(B^TX - C)$ are all in the closed left-half plane. This extends \cite[Theorem 2]{JWDSP2} to the case of input-state-output systems that are not necessarily controllable or observable. We note that, in contrast to the controllable case presented in \cite{JWDSP2}, an uncontrollable passive system will not possess a solution to the algebraic Riccati equation \eqref{eq:are} for which the eigenvalues of $A+B(D+D^T)^{-1}(B^TX - C)$ are all in the closed right-half plane, as $A$ will necessarily possess an eigenvalue in the open left-half plane corresponding to an uncontrollable mode, which will also be an eigenvalue of $A+B(D+D^T)^{-1}(B^TX - C)$.

Finally, it is also shown in \cite{pec_ba} that a system is passive if and only if it is realized as the behavior of a passive electric circuit. This extends classical network synthesis results to behaviors that need not be controllable. A constructive realization procedure is also obtained (see the sidebar titled ``\nameref{sb:pns}'').

\subsection{Necessary and sufficient conditions for reciprocity}
\label{sec:nscr}
In the paper \cite{THRSC}, a theory of reciprocal systems for behaviors corresponding to the solutions to  differential equations taking the form of \eqref{eq:gcfps} is developed, using the aforementioned representation-free definition of reciprocity. In particular, the paper \cite{THRSC} contained the following principal results:
\begin{enumerate}
\item The behavior in \eqref{eq:gcfps} constitutes a reciprocal system if and only if $PQ^T$ is symmetric.
\item If the behavior in \eqref{eq:gcfps} is reciprocal, then there always exists a subset of the variables $i$ together with a complementary subset of the variables $v$ that constitute an input to the system.
\item Any given reciprocal input-output system possesses an internally reciprocal input-state-output realization, as in \eqref{eq:ir} (for the case in which $\Sigma_e$ is the identity matrix). This extends \cite[Theorem 6]{JWDSP2} to include uncontrollable systems.
\item Any given reciprocal and passive input-output system possesses an internally reciprocal and passive input-state-output realization, as in \eqref{eq:ir}--\eqref{eq:ipir1} (again for the case in which $\Sigma_e$ is the identity matrix). This extends \cite[Theorem 7]{JWDSP2} to include uncontrollable systems.
\item The behavior corresponding to the set of solutions to equation \eqref{eq:gcfps} is passive and reciprocal if and only if it can be realized as the driving-point behavior of an electric circuit comprising resistors, inductors, capacitors, and transformers. This extends classical passive network synthesis results, as detailed in \cite{AndVong}, to include systems that are uncontrollable or not naturally formulated in terms of inputs and outputs.
\end{enumerate}
In addition to the aforementioned realization results, \cite[Theorem 8]{THISMENR} establishes bounds on the number of states of even and odd parity in internally reciprocal realizations of reciprocal input-output systems, and \cite[Theorem 9]{THISMENR} establishes bounds on the number of capacitors and inductors required to realize a given passive and reciprocal behavior. The bounds are formulated in terms of the number of uncontrollable modes (see the section titled ``\nameref{sec:ba}''), and the numbers of positive and negative eigenvalues of the Bezoutian matrix for the polynomial matrices $P$ and $Q$ in the governing differential equation \eqref{eq:gcfps}. Here, the Bezoutian matrix is the block matrix whose block entries $(Bez)_{ij}$ satisfy
\begin{equation*}
\frac{Q(z)P(w)^T - P(z)Q(w)^T}{z-w} = \sum_{i=1}^m \sum_{j=1}^m (Bez)_{ij}z^{i-1}w^{j-1},
\end{equation*}
where $m$ is the maximum degree among all entries in the matrix polynomials $P(s)$ and $Q(s)$.

Specifically, it is shown that, for a reciprocal input-state-output system, the number of states of even (respectively, odd) parity are bounded below by the number of positive (respectively, negative) eigenvalues of the Bezoutian matrix of $P$ and $Q$ plus the number of uncontrollable modes. Similarly, for a reciprocal electric circuit, the number of capacitors (respectively, inductors) are bounded below by the number of positive (respectively, negative) eigenvalues of the Bezoutian matrix of $P$ and $Q$ plus the number of uncontrollable modes. In particular, it follows that any uncontrollable mode in an internally reciprocal realization of a reciprocal input-output behavior requires at least one state of positive parity and one state of negative parity. And any uncontrollable mode in a reciprocal electric circuit requires at least one inductor and at least one capacitor.

We note that the preceding results are referring to uncontrollable modes of the driving-point behavior in the manner defined in the section titled ``\nameref{sec:ba}". However, it may still be the case that the circuit contains further \emph{internal} uncontrollable modes. See, for instance, the example in Fig.\ \ref{fig:RCcircuit} within the sidebar titled ``\nameref{sb:ecioco}". This circuit has the uncontrollable and unobservable input-state-output representation indicated in that sidebar, yet the driving-point behavior takes the form of equation \eqref{eq:gcfps} with $P(s)=2$ and $Q(s)=s+1$, whereupon there are no uncontrollable driving-point modes and the Bezoutian matrix is simply equal to $2$.

\subsection{Necessary and sufficient conditions for losslessness}
\label{sec:nscl}

It follows from the results of the preceding section that if a behavior is lossless then it is necessarily (behaviorally) controllable. 
Indeed, the behavior corresponding to the solutions to the differential equation \eqref{eq:gcfps} is lossless if and only the following algebraic conditions hold simultaneously:
\begin{enumerate}
\item $P(\lambda)Q(\bar{\lambda})^T + Q(\lambda)P(\bar{\lambda})^T \geq 0$ for all $\lambda$ in the closed right-half plane.
\item $P(i \omega)Q(-i\omega)^T + Q(i\omega)P(-i\omega)^T = 0$ for all real $\omega$.
\item The rows of $\begin{bmatrix}P(\lambda)& {-}Q(\lambda)\end{bmatrix}$ are independent for all complex $\lambda$.
\end{enumerate}
In addition, the behavior is lossless if and only if it can be realized as the driving-point behavior of an electric circuit comprising inductors, capacitors, transformers and gyrators.

\subsection{Reversible and relaxation systems}

We say that the behavior corresponding to the solutions to the differential equation \eqref{eq:gcfps} constitutes a \emph{reversible} system if it is controllable and, whenever $(i,v)$ is a trajectory from the behavior for which $i$ and $v$ are both initially at rest, then $(\hat{i},\hat{v})$ is also a trajectory from the behavior, where $\hat{i}(t) = -i(-t)$ and $\hat{v}(t) = v(-t)$ for all real $t$. It can then be shown that a dissipative system is reversible if and only if it is lossless and reciprocal or, equivalently, if and only if it is realizable as the driving-point behavior of an electric circuit comprising inductors, capacitors and transformers. In particular, the behavior corresponding to the set of solutions of the differential equation \eqref{eq:gcfps} is dissipative and reversible if and only if $PQ^T$ is symmetric and the algebraic conditions outlined in the subsection titled ``\nameref{sec:nscl}'' all hold. These results follow in a straightforward manner from the results in this article and in \cite[Section 9]{JWDSP2}. More specifically, we recall the guaranteed existence of an input $u$ comprised of a subset of the variables in $i$ together with the complementary variables in $v$, and we denote by $y$ the corresponding output. By forming matrices $\hat{P}$ and $\hat{Q}$ from the appropriate columns from the matrices $P$ and $Q$, we then obtain an alternative characterisation of the system as the set of solutions to the differential equations
\begin{equation}
\hat{P}(\tfrac{d}{dt})u = \hat{Q}(\tfrac{d}{dt})y. \label{eq:iod}
\end{equation} 
We then define $\Sigma_e$ as the signature matrix whose $j$th diagonal entry is $+1$ whenever the $j$th input corresponds to the $j$th entry in $i$, and $-1$ otherwise. It is then straightforward to show that the system in equation \eqref{eq:gcfps} is reversible (resp., reciprocal) if and only if the system in equation \eqref{eq:iod} is reversible (resp., reciprocal) with respect to the signature matrix $\Sigma_e$ in the sense of \cite{JWDSP2}.

We say that the behavior corresponding to the solutions to the differential equation \eqref{eq:gcfps} constitutes a \emph{relaxation} system if it is realizable as the driving-point behavior of an electric circuit comprising capacitors, resistors, and transformers. 

%
It then follows from the bounds on the numbers of inductors and capacitors required to realize a given behavior that were outlined in the subsection titled ``\nameref{sec:nscr}'' that a relaxation system is necessarily (behaviorally) controllable. Furthermore, the behavior corresponding to the solutions to the differential equation \eqref{eq:gcfps} constitutes a relaxation system if and only if the following algebraic conditions hold simultaneously:
\begin{enumerate}
\item $P(\lambda)Q(\bar{\lambda})^T + Q(\lambda)P(\bar{\lambda})^T \geq 0$ for all $\lambda$ in the closed right-half plane.
\item The rows of $\begin{bmatrix}P(\lambda)& {-}Q(\lambda)\end{bmatrix}$ are independent for all complex $\lambda$.
\item $PQ^T$ is symmetric and the Bezoutian matrix of $P$ and $Q$ is non-negative definite.
\end{enumerate}

A dual result to the above holds for the set of behaviors realizable as the driving-point behavior of an electric circuit comprising inductors, resistors, and transformers. The algebraic conditions are the same as for relaxation systems with the one exception that the \emph{negative} of the Bezoutian matrix of $P$ and $Q$ is non-negative definite in this case.

\section{Conclusion}
The outstanding contributions of Jan Willems' seminal Dissipative Dynamical Systems papers still prove a source of practical relevance and inspiration 50 years after their publication, and have catalysed a great body of work and contributed towards a diverse range of application areas in the interim. The papers sought to describe and generalise the notions of passivity and reciprocity, and to integrate the theories of linear optimal control and passive network analysis and synthesis. In this respect, there is very little by way of limitations of the scope and conceptual power of the papers, with the exception of the input-state-output formulation and the prevalence of controllability and observability assumptions throughout. In this article, we have highlighted such limitations, relating these to the later work of Jan Willems himself and to examples from the classes of electric circuits and mechanical systems. We have summarised recent results on passivity and reciprocity that overcome these limitations. These are based on trajectory-level and representation-free definitions of passivity and reciprocity, and lead to necessary and sufficient conditions for passivity, reciprocity, and electric circuit realizability, in terms of the differential equations describing the behavior of the system, as opposed to the system's transfer function. The trajectory-level definition of passivity extends naturally to a definition of dissipativity, and is fitting for the infinite-dimensional, time-varying, and non-linear settings. A topic for ongoing research is to extend the results presented herein into such domains.

\bibliographystyle{IEEEtran}
\bibliography{passive_ne_rec_syst_JW_persp}

\processdelayedfloats 

\sidebars 

\clearpage
\newpage

\section{Sidebar: Summary}

The dissipativity concept sits at the intersection of physics, systems theory, and control engineering, as a natural generalisation of passive systems that dissipate energy. It relates the external behavior of systems to their internal state, and connects the subjects of optimal control, algebraic Riccati equations, linear matrix inequalities, complex functions, and spectral factorization. Within control, its applications include the analysis and design of interconnected systems (such as cyber-physical systems), robustness, and the absolute stability problem, and network synthesis (of electrical, mechanical, and multi-physics systems). Dissipativity emerged as a standalone concept following the seminal Dissipative Dynamical Systems papers of Jan Willems, which drew on the work of Kalman, Youla, Anderson, Yakubovich, Popov, and others. This article details recent developments in the treatment of dissipativity and the related concept of reciprocity for systems that are not necessarily controllable and need not lend themselves naturally to an input-state-output perspective, as is the case for many physical and passive systems. We illustrate these concepts using simple electric circuit and mechanical network examples. We also draw inspiration from Jan Willems’ behavioral theory, a natural formalism for analysing physical systems that need not be controllable or characterised in terms of inputs and outputs.

\newpage
\processdelayedfloats 
\clearpage

\section[$H_{\infty}$ methods]{Sidebar: $H_{\infty}$ methods}
\label{sb:hinf}

A major contribution of the Dissipative Dynamical Systems papers of Jan Willems \cite{JWDSP1, JWDSP2} is the unification of optimal control methods, frequency domain characterisations of dissipativity, and state-space characterisations of dissipativity involving linear matrix inequalities and algebraic Riccati equations. This approach has proven conceptually powerful, in particular in the treatment of the $H_{\infty}$ optimal control problem \cite{S1}. This is concerned with the dissipativity of systems with respect to a supply rate of the form
$$
w = y^T y - \gamma^2 u^T u,
$$
for some given real number $\gamma > 0$. Indeed, the bounded-real lemma establishes results analogous to those of \cite{JWDSP2} for the case of systems that are dissipative with respect to such a supply rate. Classical treatments of the bounded-real lemma (see, for example, \cite[Section 7.2]{AndVong}) are again subject to similar controllability and observability assumptions to those detailed in this article in the context of passive systems. A treatment of the bounded-real lemma that is free from such assumptions was recently presented in \cite[Section VI]{Hughes2018}. This is for FDLTI systems in the general form
\begin{equation}
P(\tfrac{d}{dt})u = Q(\tfrac{d}{dt})y, \label{eq:gde}
\end{equation}
where $P$ and $Q$ are polynomial matrices and $Q$ is square. It is presented for the special case of $\gamma = 1$, but is easily generalised to arbitrary values of $\gamma$, thereby establishing the equivalence of the following conditions:

\begin{enumerate}
\item The system described in equation \eqref{eq:gde} is dissipative with respect to the supply rate $y^T y - \gamma^2 u^T u$. That is, for any given real number $t_0$ and any given state of the system at $t_0$, the future extracted supply
\begin{equation}
\int_{t_0}^{t_1}\gamma^2 u(t)^T u(t) - y(t)^T y(t)\mathrm{dt}
\end{equation}
is bounded above for all $t_1 \geq t_0$ and for all future trajectories of the system that share the specified state at $t_0$.
\item The following three conditions hold simultaneously
\begin{enumerate}
\item $\gamma^2 Q(\lambda) Q(\bar{\lambda})^T - P(\lambda) P(\bar{\lambda})^T \geq 0$ for all $\lambda$ in the closed right-half plane.
\item The rows of $\begin{bmatrix}P(\lambda)& {-}Q(\lambda)\end{bmatrix}$ are independent for all $\lambda$ from the closed right-half plane.
\item If $p$ is a polynomial vector and $\lambda$ a complex number that satisfy $p(s)^T(\gamma^2 Q(s)Q(-s)^T - P(s)P(-s)^{T}) = 0$ and $p(\lambda)^T\begin{bmatrix}P(\lambda)& {-}Q(\lambda)\end{bmatrix} = 0$, then $p(\lambda) = 0$.
\end{enumerate}
\item There exists an input-state-output representation of the system as in equation \eqref{eq:ssr}, and for any such input-state-output representation there exists an $X \geq 0$ such that
$$
\begin{bmatrix}-A^T X - XA - C^T C& -C^T D - XB\\ -D^T C - B^T X& \gamma^2 I - D^T D\end{bmatrix} \geq 0.
$$
\end{enumerate}

If, in addition, $\gamma^2 I - D^T D \geq 0$, then the above conditions are also equivalent to the existence of an $X \geq 0$ that satisfies the Algebraic Riccati equation:
$$
-A^T X - XA - C^T C - (C^T D + XB)(\gamma^2 I - D^T D)^{-1}(D^T C + B^T X) = 0.
$$

\newpage
\processdelayedfloats 
\clearpage

\section[Passive mechanical control]{Sidebar: Passive mechanical control}
\label{sb:mn}
The theory of dissipativity outlined in \cite{JWDSP1, JWDSP2} was influenced in part by results in electric circuit theory, and has informed the ongoing development of circuit theory through the connections with passivity. Of notable significance is the result that a FDLTI system is passive if and only if it is realizable as the driving-point behavior of an electric circuit, as described in the main body of this article. The practical significance of these results has been enhanced by the invention of the inerter \cite{mcs02}. The inerter is a one-port mechanical element with the property that the force transmitted through the element is proportional to the relative acceleration of its two terminals. It completes the force-current analogy between the electrical one-port elements resistors, inductors and capacitors (indicated in Table \ref{fig:ece}), and the mechanical one-port components dampers, springs and inerters (indicated in Table \ref{table:mechanicalcomponents}). We note that it differs from a mass as it has two terminals. In contrast, a mass has only a single terminal and is hence analogous to a grounded capacitor.

To complete the analogy between electrical and mechanical systems it is necessary also to consider mechanical analogues to transformers and gyrators (as characterised in Table \ref{fig:ece}). Mechanical analogues to transformers include gears and levers. A simple two-port realization in the form of a lever is shown in Table \ref{table:mechanicalcomponents}. It appears to be taken as implicit in the literature that multiport analogues to electric transformers are also realizable, though the authors are not aware of any proposed physical implementations of such devices. A mechanical analogue of the electric gyrator is the gyroscope.

\begin{figure}[!t]
	\centering
	\includegraphics[]{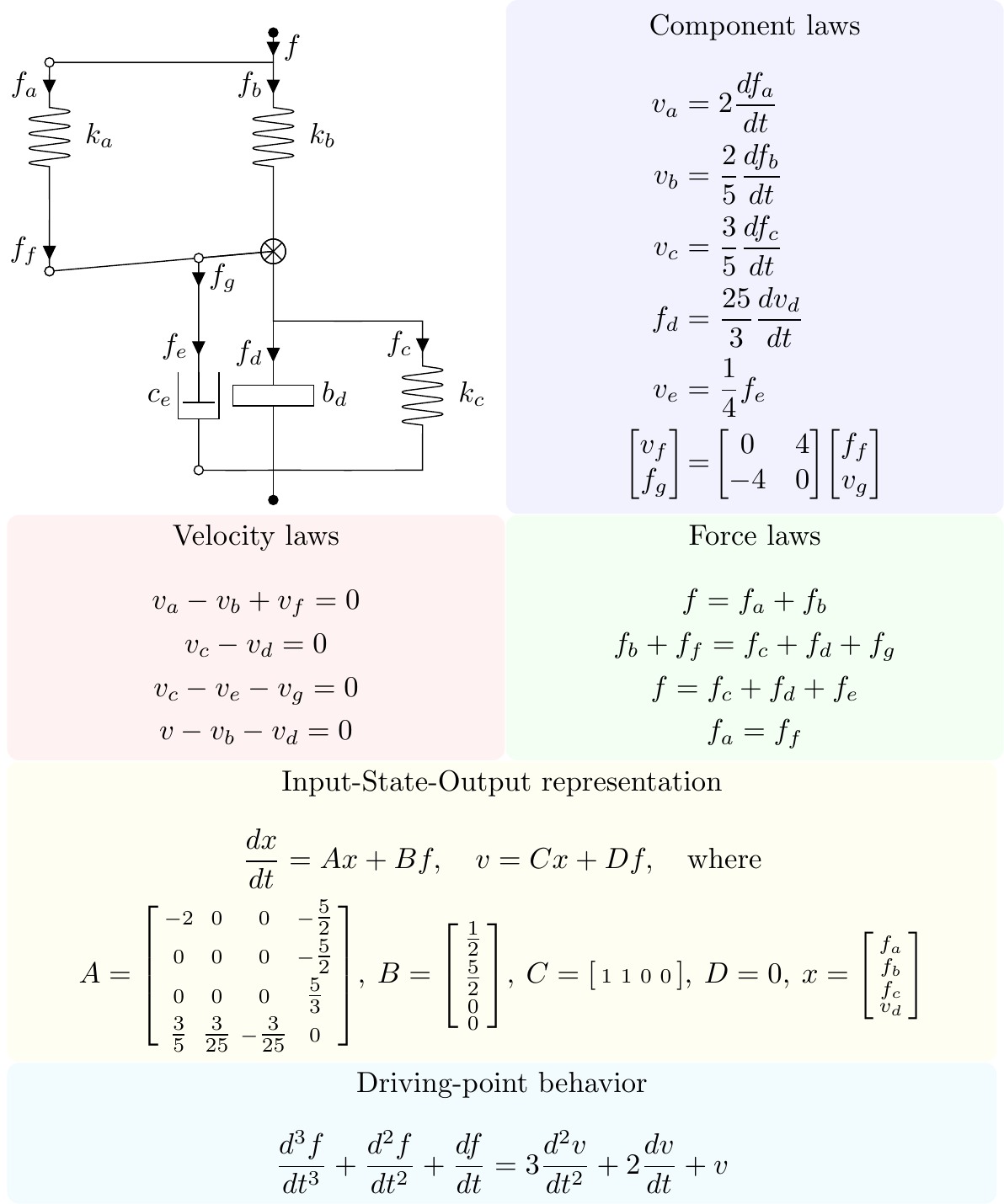}
	\caption[]{The Darlington Synthesis for the  mechanical impedance $(s^3 + s^2 + s)/(3s^2 + 2s + 1)$ (see \cite[Section IV]{mcs02} and \cite[Section 7]{pec_ba}), and the equations corresponding to the individual element laws, the velocity and force relationships, and the input-state-output and driving-point representations derived from them using the elimination theorem. This is the mechanical analogue of the same impedance as the electric circuit in Fig.\ \ref{fig:Darlington Synthesis}. Here, $f_a$ is the compressive force applied to the spring with constant $k_a$ and $v_a$ is the relative velocity of the two terminals of this spring, with the convention that $f_a v_a$ equals the instantaneous power transferred to this spring, and so forth. Similarly, $f$ is the compressive force applied to the whole network and $v$ is the relative velocity of the network terminals, with the convention that $f v$ is the instantaneous power transferred into the mechanical network.}
	\label{fig:Darlington Synthesis Mechanical}
\end{figure}

To demonstrate this force-current analogy between electrical and mechanical passive systems we include both a mechanical and an electrical realization of the same transfer function in Figures \ref{fig:Darlington Synthesis Mechanical} and \ref{fig:Darlington Synthesis} respectively. Note that the 16 equations in the 15 unknowns are equivalent according to the force-current analogy. In particular, static  equilibrium at each node (the forces on the node summing to zero) is equivalent to Kirchhoff's current law, while  compatibility (the sum of relative velocities around a loop summing to zero) is equivalent to Kirchhoff's voltage law. The elimination procedure detailed in the sidebar titled ``\nameref{sb:etape}"---for the case of the electric circuit realization---can be used to obtain the state-space representation and the driving-point behavior given in these figures, starting from either the mechanical or the electrical realization.

Taking into account the aforementioned analogy between electric circuits and mechanical networks, it follows that dissipativity theory in general, and its connections to electric circuit synthesis in particular, is of direct applicability to the design of passive mechanical controllers. Such techniques have recently found application to automotive vehicles \cite{S1}, railway vehicles \cite{S2, S3, S4}, motorcycle steering systems \cite{S5, S6}, aircraft landing systems \cite{S7}, and buildings \cite{S8, S9, S10}. We note that such mechanical applications elevate the importance of minimality and synthesis without transformers, owing to space and cost constraints. This is a research area with a number of curious results and open questions (see, for example, \cite{S11}).

\begin{table}[t]
	\centering
		\caption[]{This table shows four mechanical equivalents of electric circuit components, the damper, spring, inerter, and lever. 
		We use the convention that a positive force and a positive relative velocity correspond to compression in the case of the damper, spring and inerter. For the lever, positive $f_1$ and $f_2$ correspond to forces applied in the same direction as the associated velocities $v_1$ and $v_2$ respectively.}
	\begin{circuitikz}
		[]\draw (2,0) to[short,*-,i>=$f$] ++(0.5,0)
		to[damper] ++(1.4,0)
		to[short,-*] ++(0.5,0);
		\draw [{Bar[width=6]}-{Stealth[width=5]}] (2,-0.3)--(2.4,-0.3) node[anchor=north]{$v^+$};
		\draw [{Bar[width=6]}-{Stealth[width=5]}] (4.4,-0.3)--(4.8,-0.3) node[anchor=north]{$v^-$};
		\node at (0,0) {Damper};
		\node at (7.5,0) {$c(v^+-v^-)=f,\quad c>0$};
		\node at (0,-1.5) {Spring};
		\draw (2,-1.5) to[short,*-,i>=$f$] ++(0.5,0)
		to[spring] ++(1.4,0)
		to[short,-*] ++(0.5,0); 
		\node at (7.5,-1.5) {$k(v^+-v^-)=\frac{df}{dt},\quad L>0$};
		\draw [{Bar[width=6]}-{Stealth[width=5]}] (2,-1.8)--(2.4,-1.8) node[anchor=north]{$v^+$};
		\draw [{Bar[width=6]}-{Stealth[width=5]}] (4.4,-1.8)--(4.8,-1.8) node[anchor=north]{$v^-$};
		\node at (0,-3) {Inerter};
		\draw (2,-3) to[short,*-,i>=$f$] ++(0.5,0)
		to[inerter] ++(1.4,0)
		to[short,-*] ++(0.5,0);
		\node at (7.5,-3) {$f=b\frac{d(v^+-v^-)}{dt},\quad b>0$};
		\draw [{Bar[width=6]}-{Stealth[width=5]}] (2,-3.3)--(2.4,-3.3) node[anchor=north]{$v^+$};
		\draw [{Bar[width=6]}-{Stealth[width=5]}] (4.4,-3.3)--(4.8,-3.3) node[anchor=north]{$v^-$};
		\node at (7.5,-4.9) {$\begin{bmatrix}v_1\\f_2\end{bmatrix} = \begin{bmatrix}0&\alpha\\-\alpha&0\end{bmatrix}\begin{bmatrix}f_1\\v_2\end{bmatrix}$};
		\node at (0,-5) {Lever};
		\draw  (2,-5.3) to[short,o-] coordinate [at end] (LeverCenter) ++(2.4,0.5);
		to[short,-*] ++(0.3,0);
		\draw [short,-o] (LeverCenter)  to ($(LeverCenter)!0.7!(2,-5.3)$);
		\draw [dashed] (LeverCenter) -- ++(-2.3,0);
		\node at (LeverCenter) [forbidden sign,draw,scale=0.9]{};
		\node at (LeverCenter) [correct forbidden sign, draw,scale=0.9]{};
		\draw [{Bar[width=6]}-{Stealth[width=5]}] (2,-4.8) node[anchor=south]{$v_1$} --(2,-5.2) ;
		\draw [{Bar[width=6]}-{Stealth[width=5]}] (2.72,-4.8) node[anchor= south]{$v_2$} --(2.72,-5.1) ;
		\node at (7.5,-4.9) {$\begin{bmatrix}v_1\\f_2\end{bmatrix} = \begin{bmatrix}0&\alpha\\-\alpha&0\end{bmatrix}\begin{bmatrix}f_1\\v_2\end{bmatrix}$};
	\end{circuitikz}

	\label{table:mechanicalcomponents}
\end{table}

\newpage
\processdelayedfloats 
\clearpage

\section[Passive network analysis]{Sidebar: Passive network analysis}
\label{sb:pna}

Passive electric circuits consist of connections of linear resistors, capacitors, inductors, transformers, and gyrators, whose properties are summarised in Table \ref{fig:ece}.

The inner product of the port currents and voltages for each of the elements in Table \ref{fig:ece} is equal to the instantaneous power delivered to the element. An $n$-port RLCTG circuit then comprises an interconnection of these elements, together with a current and voltage associated with each of the $n$ ports, where the interconnection constraints are governed by Kirchhoff's current and voltage laws, and where the inner product of the driving-point currents and voltages represent the power transferred to the network. Here, Kirchhoff's current law states that the sum of currents into any vertex, or more generally, cut-set, is zero; and Kirchhoff's voltage law states that the sum of voltages around any closed loop is zero. For a circuit comprising $p$ vertices (element and external terminals) and $m = m_1 + m_2$ edges (corresponding to $m_1$ external ports and $m_2$ internal elements), then Kirchhoff's current law results in $p-1$ independent linear equations while Kirchhoff's voltage law gives $m-p+1$, and the constituent element laws result in $m_2$ differential algebraic equations. It is shown in \cite{pec_ba} that it is always the case that the internal element currents and voltages are properly eliminable (see the sidebar titled ``\nameref{sb:etape}'') to give a differential equation representation for the external behavior of the circuit comprising $m_1$ independent differential algebraic equations among the $2m_1$ driving-point currents and voltages. 

\begin{figure}[!t]
	\centering
	\includegraphics[]{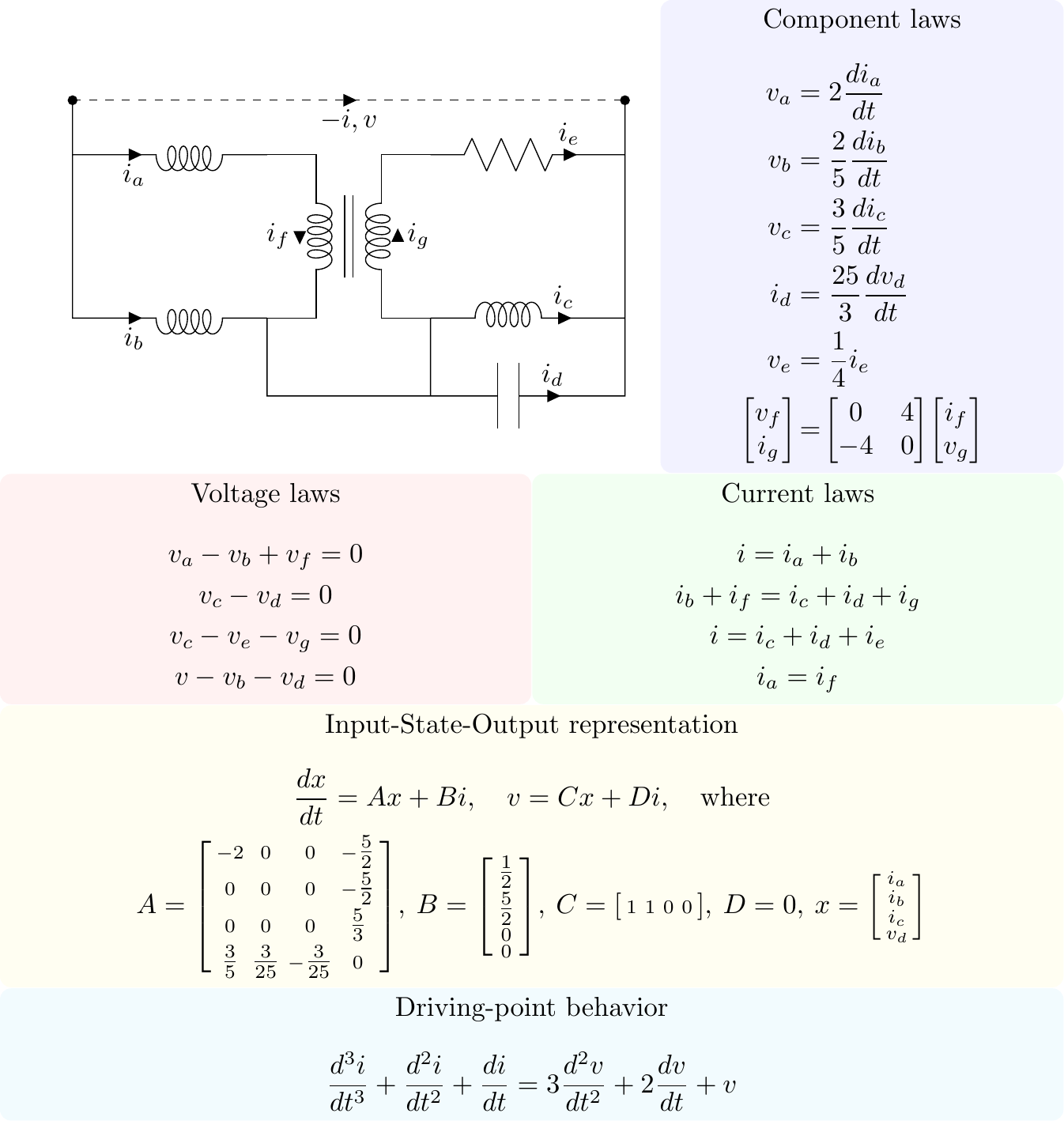}
	
	\caption[]{The Darlington Synthesis for the impedance $(s^3 + s^2 + s)/(3s^2 + 2s + 1)$ (see \cite[Section 7]{pec_ba} and \cite[Section IV]{mcs02}), and the equations corresponding to the individual element laws and Kirchhoff's current and voltage laws, and the input-state-output and driving-point representations derived from them using the elimination theorem. The sign convention is such that $v_ai_a$ is the instantaneous power transferred to the inductor whose port current is $i_a$ and similarly for $v_bi_b$ and so forth. Moreover, $vi$ is the instantaneous power transferred to the circuit.
	\newline
	\newline
	This circuit realizes the same transfer function as the mechanical network in Fig.\ \ref{fig:Darlington Synthesis Mechanical} using the force-current analogy detailed in the sidebar titled ``\nameref{sb:mn}".}
	\label{fig:Darlington Synthesis}
\end{figure}

For example, the circuit in Fig.\ \ref{fig:Darlington Synthesis} is obtained from the so-called Darlington synthesis technique (see \cite[Section IVF]{mcs02} and \cite{pec_ba}). This comprises 8 edges and 6 vertices, and Kirchhoff's laws together with the constituent element laws yield the 15 equations amongst the 16 unknowns indicated in that figure (similarly, the equivalent 15 equations in 16 unknowns in Fig.\ \ref{fig:Darlington Synthesis Mechanical} can be obtained). 
Following the elimination theorem, the resistor currents and voltages, inductor voltages, and capacitor currents can be eliminated to obtain an input-state-output representation of the circuit's external behavior, as shown in that figure. Similarly, the inductor currents and capacitor voltages can subsequently be eliminated to obtain a differential equation governing the external currents and voltages of the circuit, as also shown in that figure.

We note that Kirchhoff's current and voltage laws are power conserving, whereupon the inner product of the driving-point current and voltages represents the power transferred to the network, which is the sum of the powers transferred to each of the elements within the network. However, it is sometimes the case that this interconnection violates the well-posedness assumption used in \cite[Section 4]{JWDSP2} to establish the dissipativity of interconnections of dissipative systems, which requires that the state of the interconnected system comprises the union of the states of the subsystems and that the inputs to the subsystems be uniquely determined by the input to the overall system and its initial state. For example, consider the circuit comprising an interconnection of two identical three-port transformers with turns ratio $T = \begin{bmatrix}1& 1\end{bmatrix}$, whose first and second ports are both connected together. This results in the element laws
\begin{equation*}
\begin{bmatrix}\hat{v}_1\\ \hat{v}_2\\ i_1\end{bmatrix} = \begin{bmatrix}0& 0& 1\\ 0& 0& 1\\ -1& -1& 0\end{bmatrix}\begin{bmatrix}\hat{i}_1\\ \hat{i}_2\\ v_1\end{bmatrix}, \text{ and } \begin{bmatrix}\hat{v}_3\\ \hat{v}_4\\ i_2\end{bmatrix} = \begin{bmatrix}0& 0& 1\\ 0& 0& 1\\ -1& -1& 0\end{bmatrix}\begin{bmatrix}\hat{i}_3\\ \hat{i}_4\\ v_2\end{bmatrix},
\end{equation*}
together with the interconnection laws $\hat{i}_3 = -\hat{i}_1$, $\hat{i}_4 = -\hat{i}_2$, $\hat{v}_3 = \hat{v}_1$, and $\hat{v}_4 = \hat{v}_2$. The internal currents and voltages $\hat{i}_j$ and $\hat{v}_j$ ($j = 1,\ldots , 4$) can then be eliminated to obtain the external behavior $i_1 = -i_2$ and $v_2 = v_1$, corresponding to a passive system whose external behavior is equivalent to a two-port transformer with turns ratio $1$. However, the corresponding internal currents and voltages can be shown to satisfy $\hat{i}_2 = i_2 - \hat{i}_1$, $\hat{i}_3 = -\hat{i}_1$, $\hat{i}_4 = \hat{i}_1 - i_2$, and $\hat{v}_1 = \hat{v}_2 = \hat{v}_3 = \hat{v}_4 = v_1$, where $\hat{i}_1$ itself can take on any value independently of the external currents and voltages $i_1, i_2, v_1$ and $v_2$.

\begin{table}[t]
	\centering
		\caption[]{Electric circuit components. The behaviors of the elements correspond to the set of locally square integrable $(i,v)$ that satisfy the differential equations in the right-hand column (in a weak sense). Each of the elements is passive in accordance with the definition in the section titled ``\nameref{sec:passivity}'', as is any interconnection of these elements. Moreover, the resistor, inductor, capacitor, and transformer elements are all reciprocal in accordance with the definition in the section titled ``\nameref{sec:recip}'', as is any interconnection of these elements.
		\newline
		\newline
		\quad In general the transformer contains $n\geq 2$ ports. The displayed image is typically used for the two port transformer. Here, $v_1,v_2$ is a partition (compatible with $T$) of the port voltages, and $i_1,i_2$ is a compatible partition of the port currents.}
	\begin{circuitikz}
		[]\draw (2,0) to[short,*-,l_=$v^+$] ++(0.3,0)
		to[R,i=$i$] ++(2.3,0)
		to[short,-*,l_=$v^-$] ++(0.3,0);
		\node at (0,0) {Resistor};
		\node at (7.5,0) {$v^+-v^-=iR,\quad R>0$};
		\node at (0,-1.5) {Inductor};
		\draw (2,-1.5) to[short,*-,l_=$v^+$] ++(0.3,0)
		to[L,i=$i$] ++(2.3,0)
		to[short,-*,l_=$v^-$] ++(0.3,0); 
		\node at (7.5,-1.5) {$v^+-v^-=L\frac{di}{dt},\quad L>0$};
		\node at (0,-3) {Capacitor};
		\draw (2,-3) to[short,*-,l_=$v^+$] ++(0.3,0)
		to[C,i=$i$] ++(2.3,0)
		to[short,-*,l_=$v^-$] ++(0.3,0);
		\node at (7.5,-3) {$i=C\frac{d(v^+-v^-)}{dt},\quad C>0$};
		\ctikzset{voltage/bump b=0.2}
		\node at (0,-5.3) {Transformer};
		\draw (2.2,-5.3) node[transformer core,right] (T) {} 
		(T.A1) to[short,-*] ++(-0.2,0)
		(T.A2) to[short,-*] ++(-0.2,0)
		(T.B1) to[short,-*] ++(0.2,0)
		(T.B2) to[short,-*] ++(0.2,0)
		(T.A2) [open,v^=$v_1$] to (T.A1)
		(T.B2) [open,v_=$v_2$] to (T.B1);
		\node [currarrow,rotate=-90,label = left:$i_1$] at ($(T.outer dot A1)!0.5!(T.outer dot A2)$){};
		\node [currarrow,rotate=-90,label = right:$i_2$] at ($(T.outer dot B1)!0.5!(T.outer dot B2)$){};
		\node at (7.5,-5.2) {$\begin{bmatrix}v_1\\i_2\end{bmatrix} = \begin{bmatrix}0&T^T\\-T&0\end{bmatrix}\!\begin{bmatrix}i_1\\v_2\end{bmatrix}$};
		\node at (7.5,-6) {\small$T\in\mathbb{R}^{m\times n}$};
		\draw (2.2,-8.3) node[gyrator,right] (G) {}
		(G.A1) to[short,-*] ++(-0.2,0)
		(G.A2) to[short,-*] ++(-0.2,0)
		(G.B1) to[short,-*] ++(0.2,0)
		(G.B2) to[short,-*] ++(0.2,0)
		(G.A2) [open,v^=$v_1$] to (G.A1)
		(G.B2) [open,v_=$v_2$] to (G.B1);
		\node [currarrow,rotate=-90,label = left:$i_1$] at ($(G.outer dot A1)!0.5!(G.outer dot A2)$){};
		\node [currarrow,rotate=-90,label = right:$i_2$] at ($(G.outer dot B1)!0.5!(G.outer dot B2)$){};
		\node at (0,-8.3) {Gyrator};
		\node at (7.5,-8.3) {$\begin{bmatrix}v_1\\v_2\end{bmatrix} = \begin{bmatrix}0&-1\\1&0\end{bmatrix}\begin{bmatrix}i_1\\i_2\end{bmatrix}$};
	\end{circuitikz}
	\label{fig:ece}
\end{table}

\newpage
\processdelayedfloats 
\clearpage

\clearpage
\section[The elimination theorem and proper elimination]{Sidebar: The elimination theorem and proper elimination}
\label{sb:etape}

A key concept in control theory in general, and dissipativity theory in particular, is the interplay between the internal and external properties of a given system. For example, the concept of state and input-state-output representations of an input-output system; the relationship between the driving-point currents and voltages of an electric circuit, and its internal element currents and voltages; and the notions of internally passive and internally reciprocal realizations of passive and reciprocal systems. To pass from the internal behavior of a given system to its external behavior requires elimination of internal variables, a subject formalised in the elimination theorem of Jan Willems' behavioral theory. This is relevant to, for example, obtaining the differential equation describing the driving-point behavior of a given electric circuit from the Kirchhoff's laws and constituent element laws governing its internal behavior. Or obtaining the input-output behavior of a given input-state-output system. 

A key consequence of the elimination theorem is that whenever the internal behavior is characterised by a FDLTI differential equation, then the external behavior will also by characterised by a FDLTI differential equation, with the possible addition of some smoothness constraints when considering locally square integrable signals as the solution space. The concept of proper elimination, as studied in \cite{S1}, is relevant to determining those situations in which no such smoothness constraints arise.

Specifically, for the FDLTI system $R_1(\tfrac{d}{dt})w = R_2(\tfrac{d}{dt})l$, then the variables $w$ are constrained to satisfy the differential equation $R(\tfrac{d}{dt})w = 0$, where $R(s) = M(s)R_1(s)$ for a polynomial matrix $M(s)$ whose rows form a basis for the left syzygy of $R_2(s)$. Here, a left syzygy is the polynomial matrix generalisation of the left null space of a real matrix. Specifically, $M(s)$ is a polynomial matrix for which the rows of $M(\lambda)$ are independent for all complex values of $\lambda$, the number of rows of $M(s)$ is equal to the maximal dimension of the left null space of $R_2(\lambda)$ over all complex values of $\lambda$, and $M(s)R_2(s) = 0$. In fact, the set of all infinitely differentiable $w$ such that there exists an infinitely differentiable $l$ satisfying the originating differential equation is equal to the set of all infinitely differentiable $w$ that satisfy $R(\tfrac{d}{dt})w = 0$. If, in addition, the variables $l$ are \emph{properly eliminable} (see \cite{S1}), then it is also the case that the set of all \emph{locally square integrable} $w$ such that there exists a \emph{locally square integrable} $l$ satisfying the originating differential equation is equal to the set of all \emph{locally square integrable} $w$ that satisfy $R(\tfrac{d}{dt})w = 0$.

We note that the state of a given input-state-output system as in equation \eqref{eq:ssr} is always properly eliminable. Indeed, the state is always uniquely determined by its initial value $x(0)$ and the input as the absolutely continuous function that satisfies the variation of the constants formula. Moreover, it is shown in \cite{pec_ba} that the internal currents and voltages of any given electric circuit are always properly eliminable. However, for some FDLTI systems, there are subsets of variables which are not properly eliminable. An example is provided in \cite{S1} of the behavior described by the differential equations $w_1 - w_2 = 0$ and $\tfrac{dw_1}{dt} = l$, for which the external behavior corresponds to the set of solutions to the equation $w_1 - w_2 = 0$ for which the derivative of $w_1$ is locally integrable. In this case, the requirement that $l$ be locally integrable precludes the situation in which $w_1$ is a step function.

To illustrate the elimination theorem we consider the electric circuit of the Darlington synthesis given in Fig.\ \ref{fig:Darlington Synthesis}, though it is equally applicable to its mechanical analogue in Fig.\ \ref{fig:Darlington Synthesis Mechanical}. The 15 equations may be written in the form $R_1(\tfrac{d}{dt})w = R_2(\tfrac{d}{dt})l$, where
\begin{equation*}
R_1(s){=}\left[\begin{smallmatrix}2\,s&0&0&0&0&0\\ \noalign{\medskip}0&{
	\frac {2\,s}{5}}&0&0&0&0\\ \noalign{\medskip}0&0&{\frac {3\,s}{5}}&0&0
&0\\ \noalign{\medskip}0&0&0&{\frac {25\,s}{3}}&0&0
\\ \noalign{\medskip}0&0&0&0&0&0\\ \noalign{\medskip}0&0&0&0&0&0
\\ \noalign{\medskip}0&0&0&0&0&0\\ \noalign{\medskip}0&0&0&0&0&0
\\ \noalign{\medskip}0&0&0&1&0&0\\ \noalign{\medskip}0&0&0&0&0&0
\\ \noalign{\medskip}0&0&0&-1&0&1\\ \noalign{\medskip}1&1&0&0&-1&0
\\ \noalign{\medskip}0&1&-1&0&0&0\\ \noalign{\medskip}0&0&1&0&-1&0
\\ \noalign{\medskip}-1&0&0&0&0&0
\end{smallmatrix}\right]\! ,
R_2(s){=}\left[\begin{smallmatrix}1&0&0&0&0&0&0&0&0&0
\\ \noalign{\medskip}0&1&0&0&0&0&0&0&0&0\\ \noalign{\medskip}0&0&1&0&0
&0&0&0&0&0\\ \noalign{\medskip}0&0&0&1&0&0&0&0&0&0
\\ \noalign{\medskip}0&0&0&0&-{\frac{1}{4}}&1&0&0&0&0
\\ \noalign{\medskip}0&0&0&0&0&0&0&0&1&-4\\ \noalign{\medskip}0&0&0&0&0
&0&4&1&0&0\\ \noalign{\medskip}1&-1&0&0&0&0&0&0&1&0
\\ \noalign{\medskip}0&0&1&0&0&0&0&0&0&0\\ \noalign{\medskip}0&0&-1&0&0
&1&0&0&0&1\\ \noalign{\medskip}0&1&0&0&0&0&0&0&0&0
\\ \noalign{\medskip}0&0&0&0&0&0&0&0&0&0\\ \noalign{\medskip}0&0&0&1&0
&0&-1&1&0&0\\ \noalign{\medskip}0&0&0&-1&-1&0&0&0&0&0
\\ \noalign{\medskip}0&0&0&0&0&0&-1&0&0&0
\end{smallmatrix}\right]\! ,
w=\left[\begin{smallmatrix}
i_a\\\noalign{\medskip}i_b\\\noalign{\medskip}i_c\\\noalign{\medskip}v_d\\\noalign{\medskip}i\\\noalign{\medskip}v
\end{smallmatrix}\right]\! , l=\left[\begin{smallmatrix}v_a\\\noalign{\medskip}v_b\\\noalign{\medskip}v_c\\\noalign{\medskip}i_d\\\noalign{\medskip}i_e\\
\noalign{\medskip}v_e\\\noalign{\medskip}i_f\\\noalign{\medskip}i_g\\\noalign{\medskip}v_f\\\noalign{\medskip}v_g\end{smallmatrix}\right]\! .
\end{equation*}
In this case, we find that the rows of the matrix
\begin{equation*}
M(s)=\left[\begin{smallmatrix}{\frac{1}{2}}&0&0&0&-2&{\frac{
		1}{2}}&-{\frac{1}{2}}&-{\frac{1}{2}}&2&2&-{\frac{1}{2}}&-{\frac{1}{2}}
&{\frac{1}{2}}&{\frac{1}{2}}&-{\frac{5}{2}}\\ \noalign{\medskip}0&{
	\frac{5}{2}}&0&0&0&0&0&0&0&0&-{\frac{5}{2}}&0&0&0&0
\\ \noalign{\medskip}0&0&{\frac{5}{3}}&0&0&0&0&0&-{\frac{5}{3}}&0&0&0&0
&0&0\\ \noalign{\medskip}0&0&0&{\frac{3}{25}}&0&0&{\frac{3}{25}}&0&0&0
&0&0&-{\frac{3}{25}}&0&{\frac{3}{5}}\\ \noalign{\medskip}0&0&0&0&0&0&0
&0&0&0&0&-1&0&0&0
\end{smallmatrix}\right],
\end{equation*}
form a basis for the left syzygy of $R_2(s)$, and we let $R(s)=M(s)R_1(s)$. Note that $M(s)$ is a real matrix but in general will be polynomial. Thus, following the elimination theorem, we can eliminate the inductor voltages, capacitor currents, and the currents and voltages associated with the resistors and transformers to obtain the relationship $R(\tfrac{d}{dt})w=0$ for the circuit of Fig.\ \ref{fig:Darlington Synthesis}, which yields the input-state-output representation given in that figure. Furthermore, by writing the input-state-output representation in the form  $R_1\!\left(\tfrac{d}{dt}\right)w=R_2\!\left(\tfrac{d}{dt}\right)l$, where 
\begin{equation*}
R_1(s)=\begin{bmatrix}
0&-\frac{1}{2}\\
0 &-\frac{5}{2}\\
0&0\\
0&0\\
1&0
\end{bmatrix},R_2(s)=\begin{bmatrix}
{-}s-2 & 0 & 0 & -\frac{5}{2}\\
0 & -s & 0 & -\frac{5}{2}\\
0 & 0 & -s & \frac{5}{3}\\
\frac{3}{5} & \frac{3}{25} &-\frac{3}{25} &-s\\
1 & 1 & 0 & 0
\end{bmatrix},w=\begin{bmatrix}i\\v\end{bmatrix},l=\begin{bmatrix}i_a\\i_b\\i_c\\v_d\end{bmatrix},
\end{equation*}
it is possible to eliminate the states by again finding the matrix $M(s)$ whose rows form a basis for the left syzygy of $R_2(s)$. In this case
\begin{equation*}
M(s) = \begin{bmatrix}{s}^{2}-s&{s}^{2}+s+{\frac{2}{5}}&{\frac
	{3}{5}}&-5\,s&{s}^{3}+{s}^{2}+s
\end{bmatrix},
\end{equation*}
and we let $R(s) = M(s)R_1(s)$, whereupon $R\left(\tfrac{d}{dt}\right)w = 0$ obtains the driving-point behavior of the electric circuit shown in Fig.\ \ref{fig:Darlington Synthesis}. The same procedure also yields the driving-point behavior of the mechanical network in Fig.\ \ref{fig:Darlington Synthesis Mechanical}.

\newpage
\processdelayedfloats 
\clearpage

\clearpage
\section[Controllability and observability of dissipative and cyclo-dissipative systems]{Sidebar: Controllability and observability of dissipative and cyclo-dissipative systems}
\label{sb:cdac}

Alongside the concept of dissipativity sits the concept of cyclo-dissipativity (see, for example, \cite{S1, S2} and the references therein). This is typically characterised by systems that possess any storage function (not necessarily non-negative). In other words, a system is said to be cyclo-dissipative with respect to a supply rate $w$ if there exists any function $S$ of the state of the system that satisfies the dissipation inequality in equation \eqref{eq:di}. Somewhat confusingly, the cyclo-dissipativity concept is referred to as dissipativity by certain authors, as exemplified by many of the references in this sidebar.

There have been a number of notable papers on the cyclo-dissipativity of systems that are not necessarily controllable. For example, \cite{S3} and \cite{S4} build on earlier results relating to dissipative systems from \cite{S5}, which extended the positive-real lemma for input-state-output systems, as in equation \eqref{eq:ssr}, to include uncontrollable systems, in the case that $A$ has all of its eigenvalues in the open left-half plane. And \cite{S6, S7, S8} approach this problem from a behavioral perspective, building on the formalism of quadratic differential forms introduced in \cite{S9}. 

However, we note that questions still remain concerning cyclo-dissipativity of uncontrollable systems, for example regarding the observability of the storage function \cite{S10}. Indeed, the formulation of cyclo-dissipativity in an input-state-output setting is representation dependent. This is illustrated by the following example adapted from \cite{S10}. Consider the input-state-output system in the standard form of equation \eqref{eq:ssr} with $A = C = 1$ and $B = D = 0$ and the supply rate $w(u,y) = uy$. Then there doesn't exist a storage function $S$ of the state $x$ that satisfies the dissipation inequality in equation \eqref{eq:di}. However, the same input-output behavior is also exhibited by the input-state-output system, as in equation \eqref{eq:ssr}, with
$$
A = \begin{bmatrix}1& 0\\ 0& -1\end{bmatrix}, \hspace{0.3cm} B = \begin{bmatrix}0\\ 1\end{bmatrix}, \hspace{0.3cm} C = \begin{bmatrix}1& 0\end{bmatrix} \text{ and } D = 0,
$$
and for this system, $\int_{t_0}^{t_1}u(t)y(t) \mathrm{dt} = S(x(t_1)) - S(x(t_0))$ for the storage function $S(x) = x_1 x_2$. In contrast, we note that a consequence of the results presented in this paper is that the definition of passivity is not representation dependent.

\newpage
\processdelayedfloats 
\clearpage

\clearpage
\section[Controllability and observability of electric circuits]{Sidebar: Controllability and observability of electric circuits}
\label{sb:ecioco}

A key message of this article is that conventionally held assumptions of controllability, observability, and even stabilisability of passive input-state-output behaviors are not generally applicable to electric circuits. Indeed, there are a number of significant circuits that arise in classical network synthesis techniques whose input-state-output behaviors violate these assumptions. One such example is the network in Fig. \ref{fig:Darlington Synthesis} obtained from the so-called Darlington synthesis method. Taking the states of this network to be the inductor currents and capacitor voltages, then the input-state-output representation for this behavior (as described in Fig.\ \ref{fig:Darlington Synthesis}) is neither controllable nor observable. On the other hand, we note that the driving-point behavior of this network is in fact (behaviorally) controllable.

Other famous examples include the 
Bott-Duffin networks and their simplifications (see \cite{HugSmSP, HugNa, HugGMF}). In fact, it is shown in \cite{HugNa} that the driving-point behavior of the Bott-Duffin network is not (behaviorally) controllable. Moreover, taking the states of this network to be the capacitor voltages and inductor currents, then the input-state-output representation of its behavior is neither controllable nor observable. The simplified Bott-Duffin networks identified by Reza, Pantell, Fialkow, and Gerst (see \cite{HugGMF}) possess similar properties as concerns controllability and observability. Nevertheless, it is shown in \cite{HugSmSP} that the Bott-Duffin networks actually contain the least possible number of inductors and capacitors among the class of \emph{series-parallel} networks for the realization of certain positive-real functions. Similarly, the simplified networks of Reza, Pantell, Fialkow, and Gerst contain the least possible number of inductors and capacitors among the entire class of resistor-inductor-capacitor networks for the realization of certain positive-real functions.

\begin{figure}[t]
	\centering
	\normalsize
	\includegraphics[]{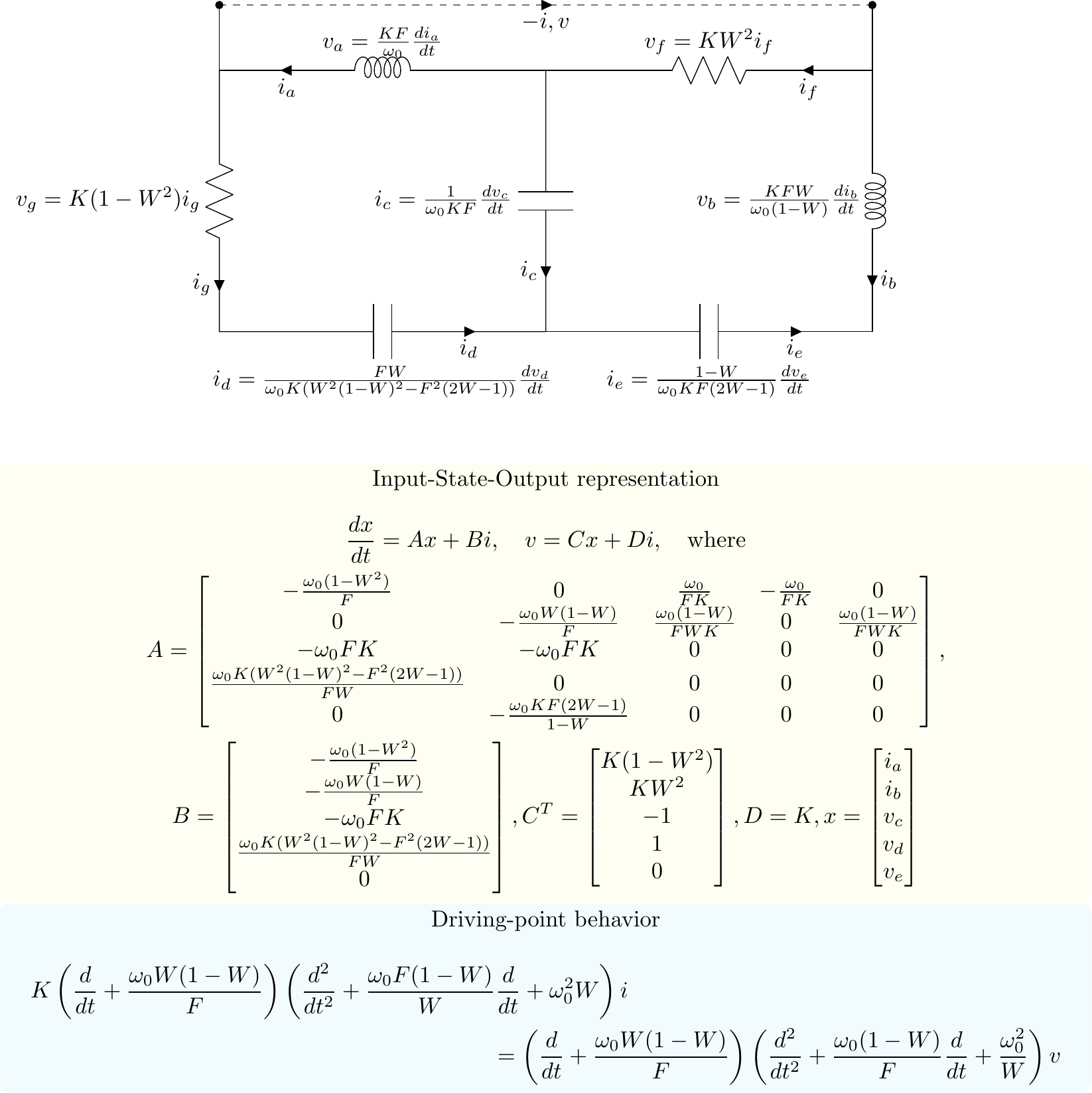}
		\caption[]{A circuit realizing a general biquadratic minimal function, adapted from \cite{HugGMF}. The component equations are labelled, and the parameters satisfy $K,\omega_0 >0$, $0<F<W(1-W)/\sqrt{2W-1}$, and $1/2<W<1$. The sign convention is such that $v_ai_a$ is the instantaneous power transferred to the inductor whose port current is $i_a$, and similarly for $v_bi_b$ and so forth, while $vi$ is the instantaneous power transferred to the circuit. Alongside the circuit diagram is an input-state-output representation and the driving-point behavior. The state-space representation is unobservable and is not stabilisable, while the driving-point behavior is (behaviorally) uncontrollable, yet is stabilisable.} 
	\label{fig:RLCHugGMF}
\end{figure}

The paper \cite{HugGMF} presented other circuits that, similar to those identified by Reza, Pantell, Fialkow, and Gerst, also contain the least possible number of inductors and capacitors for the realization of their driving-point impedances among the entire class of resistor-inductor-capacitor circuits. One example is shown in Fig.\ \ref{fig:RLCHugGMF} which has been adapted from the right-hand side of \cite[Fig.\ 3]{HugGMF}. We note a sign error in two entries of the second row of $A$ in Section II of that paper, which we correct here. We note that this input-state-output representation is not observable, since $\tilde{x} = \left[\begin{smallmatrix}{-}(1-W) W^2\:&(1-W)^2 W\:&{-}KF^2(2W-1)\:&K(W^2(1-W)^2{-}F^2(2W-1))\:&KF^2(2W-1)\end{smallmatrix}\right]^T$ satisfies $A\tilde{x} = {-}\frac{\omega_0W(1-W)}{F}\tilde{x}$ and $C\tilde{x}=0$. Moreover, $\hat{x}=\left[\begin{smallmatrix}0&0&(2W-1)(W^2(1-W)^2{-}F^2(2W-1))\;&F^2W(2W-1)\;&{-}(1-W)(W^2(1-W)^2{-}F^2(2W-1))\end{smallmatrix}\right]^T$ satisfies $\hat{x}^T A = 0$ and $\hat{x}^TB = 0$, so this input-state-output representation is not controllable, and in fact it is not even stabilisable. 

The driving-point behavior of this circuit can be obtained by the elimination procedure described in the sidebar titled ``\nameref{sb:etape}'' and in an entirely analogous manner to that described in the Darlington synthesis example in that sidebar. This driving-point behavior is given in Fig\ \ref{fig:RLCHugGMF}, and we note it is (behaviorally) uncontrollable, yet is (behaviorally) stabilisable, despite the fact that the input-state-output representation is not stabilisable.

Also of interest from \cite{HugGMF} is the circuit on the left-hand side of \cite[Fig.\ 3]{HugGMF}, which, in contrast to the examples in this article, does \emph{not} possess an input-state-output representation in which the states correspond to the set of inductor currents and capacitor voltages. This is because that circuit contains a loop comprising only capacitors, hence the capacitor voltages are linearly dependent.

\begin{figure}
	\centering
	\includegraphics[]{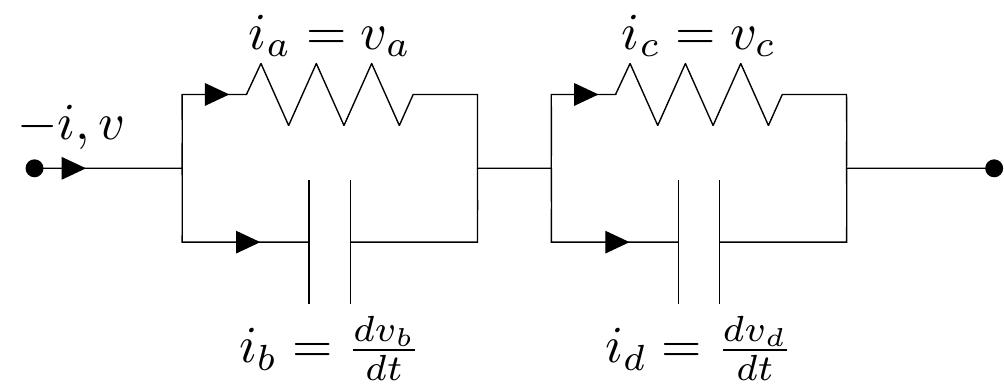}
	\caption[]{An example of a circuit that is behaviorally controllable but with uncontrollable internal modes.}
	\label{fig:RCcircuit}
\end{figure}

We note our focus is on driving-point controllability and observability of electric circuits, which may have further \emph{internal} uncontrollable modes. This is illustrated by the example in Fig.\ \ref{fig:RCcircuit}. It is straightforward to show that this circuit possesses the input-state-output representation

\begin{equation*}
\begin{gathered}
\begin{bmatrix}\tfrac{dv_b}{dt}\\ \tfrac{dv_d}{dt}\end{bmatrix} = \begin{bmatrix}-1 & 0\\ 0 & -1 \end{bmatrix}\begin{bmatrix}v_b\\v_d\end{bmatrix} + \begin{bmatrix}1 \\ 1 \end{bmatrix}i\\
v = \begin{bmatrix}1 & 1\end{bmatrix}\begin{bmatrix}v_b \\ v_d \end{bmatrix},
\end{gathered}
\end{equation*}
which is neither observable or controllable. However, by eliminating the states $v_b$ and $v_d$ in a manner similar to the previous examples, we find that the driving-point behavior is given by
\begin{equation*}
\frac{dv}{dt} + v = 2i,
\end{equation*}
which is (behaviorally) controllable. 

\tikzset{ RLCfig/.pic={
		\ctikzset{bipoles/thickness=1}
		\coordinate (topL) at (0,0);
		\draw (topL) to[short,-*] coordinate [at end] (TOPL) ++(0,1);
		\draw (topL) to[L,i<_=$i_a$,l=${v_a = \tfrac{KF}{\omega_0}\tfrac{di_a}{dt}}$] coordinate [at end] (topC) ++(5,0)
		to[R,i_<=$i_f$,l=${v_f= KW^2i_f}$] coordinate [at end] (topR) ++(5,0)
		to[L,i=$i_b$,l_=${v_b = \tfrac{KFW}{\omega_0(1-W)}\tfrac{di_b}{dt}}$] ++(0,-4)
		to[C,i<=$i_e$,l=${\.i_e=\tfrac{1-W}{\omega_0KF(2W-1)}\tfrac{dv_e}{dt}}$] coordinate [at end] (bottomC) ++(-5,0);
		\draw (bottomC) to[C,i<=$i_c$,l=${i_c=\tfrac{1}{\omega_0 KF}\tfrac{dv_c}{dt}}$] (topC)
		(bottomC) to[C,i<=$i_d$,l=${i_d=\tfrac{FW}{\omega_0 K(W^2(1-W)^2 - F^2(2W-1))}\tfrac{dv_d}{dt}}$] ++(-5,0)
		to[R,i<=$i_g$,l=${v_g=K(1-W^2) i_g}$]  (topL);
		\draw (topR) to[short,-*] coordinate [at end] (TOPR) ++(0,1);
		\draw [dashed] (TOPL)--(TOPR);
		\node [currarrow,label=below:${-i,v}$] at ($(TOPL)!0.5!(TOPR)$) {};
}}

\newpage
\processdelayedfloats 
\clearpage

\clearpage
\section[Passive network synthesis]{Sidebar: Passive network synthesis}
\label{sb:pns}
The results detailed in this article lead naturally to a passive network synthesis method that extends the reactance extraction method of \cite{YoulaTissi, AndVong} to those behaviors that are uncontrollable. Specifically, this provides an explicit realization for any given passive behavior as the driving-point behavior of a circuit comprising resistors, inductors, capacitors, transformers, and gyrators, and an explicit realization for any given passive and reciprocal behavior as the driving-point behavior of a circuit comprising resistors, inductors, capacitors, and transformers. For the case of a passive (resp., passive and reciprocal) input-output system, this synthesis method follows from the existence of an internally passive (resp., internally passive and reciprocal) realization, which can be realized using the methods outlined in \cite[Section VI(2)]{S1}. For those cases in which $i$ does not have the property of an input, an explicit electric circuit realization is obtained, taking the form of \cite[Fig.\ 2]{THRSC}, as explained in the proof of Theorem 17 of that paper.

\newpage
\section{Acknowledgements}
We would like to thank the anonymous reviewers for their very detailed and useful comments that have helped us improve this article to a great extent. In particular we thank the anonymous reviewer 1 for suggesting the simple example we now include in Fig.\ \ref{fig:RCcircuit}.

\newpage
\section{Author Biography}
Timothy H.\ Hughes received the M.Eng.\ degree in Mechanical Engineering, and the Ph.D.\ degree in Control Engineering, from the University of Cambridge, U.K., in 2007 and 2014, respectively. From 2007 to 2010 he was employed as a mechanical engineer at The Technology Partnership, Hertfordshire, U.K; and from 2013 to 2017 he was a Henslow Research Fellow at the University of Cambridge, funded by the Cambridge Philosophical Society. He is currently a lecturer in the Department of Mathematics at the University of Exeter, U.K.

Edward H.\ Branford is a PhD student in the Department of Mathematics at the University of Exeter, U.K., where he is supervised by Timothy H.\ Hughes and Stuart B.\ Townley. He received the B.A.\ degree in Physics from the University of Oxford, U.K., in 2016. He received the M.Sc.\ degree in Systems, Control and Signal Processing from the University of Southampton, U.K., in 2019.
\end{document}